\documentclass[11pt]{article}
\usepackage{amssymb}
\usepackage{amsthm}
\usepackage{amsmath,amsfonts,amssymb}

\usepackage[mathscr]{eucal}
\usepackage{exscale}
\usepackage{natbib}
\usepackage{bm}
\usepackage{eqlist}
\usepackage[dvipsnames]{color}
\usepackage[Lenny]{fncychap}
\usepackage{multirow}
\usepackage{graphicx}
\usepackage{algpseudocode}
\usepackage{algorithm}
\usepackage{framed}

\usepackage{tabu}

\newtheorem{theorem}{Theorem}[section]

\newtheorem{definition}[theorem]{Definition}
\newtheorem{remark}[theorem]{Remark}

\font\bigbf=cmbx10 scaled \magstep3

\begin{document}

\title{\bigbf Computing dynamic user equilibria on large-scale networks: From theory to software implementation}

\author{Ke Han$^{a}\thanks{Corresponding author; e-mail: k.han@imperial.ac.uk}$ 
\quad Gabriel Eve$^{a}\thanks{e-mail: gabrielfeve@gmail.com}$ 
\quad Terry L. Friesz$^{b}\thanks{e-mail: tfriesz@psu.edu}$ 
\\\\
$^{a}$\textit{Department of Civil and Environmental Engineering}\\
\textit{Imperial College London, SW7 2BU, UK}\\ 
$^{b}$\textit{Department of Industrial and Manufacturing Engineering}\\
\textit{Pennsylvania State University, PA 16803, USA}\\
}
\date{}

\maketitle

\begin{abstract}
Dynamic user equilibrium (DUE) is the most widely studied form of dynamic traffic assignment, in which road travelers engage in a non-cooperative Nash-like game with departure time and route choices. DUE models describe and predict the time-varying traffic flows on a network consistent with traffic flow theory and travel behavior. This paper documents theoretical and numerical advances in synthesizing traffic flow theory and DUE modeling, by presenting a holistic computational theory of DUE with numerical implementation encapsulated in a MATLAB software package. In particular, the dynamic network loading (DNL) sub-problem is formulated as a system of differential algebraic equations based on the fluid dynamic model, which captures the formation, propagation and dissipation of physical queues as well as vehicle spillback on networks. Then, the fixed-point algorithm is employed to solve the DUE problems on several large-scale networks. We make openly available the MATLAB package, which can be used to solve DUE problems on user-defined networks, aiming to not only help DTA modelers with benchmarking a wide range of DUE algorithms and solutions, but also offer researchers a platform to further develop their own models and applications. Updates of the package and computational examples are available at \textcolor{blue}{https://github.com/DrKeHan/DTA}.
\end{abstract}


\section{Introduction}\label{secIntro}

This paper is concerned with a class of models known as {\it dynamic user equilibrium} (DUE). DUE problems have been studied within the broader context of {\it dynamic traffic assignment} (DTA), which is viewed as the modeling of time-varying flows on traffic networks consistent with established travel demand and traffic flow theory.

DTA models, from the early 1990s onward, have been greatly influenced by Wardrop's principles \citep{Wardrop1952}. 
\begin{itemize}
\item Wardrop's first principle, also known as the {\it user optimal} principle, views travelers as Nash agents competing on a network for road capacity. Specifically, the travelers selfishly seek to minimize their own travel costs by adjusting route choices. A user equilibrium is envisaged where the travel costs of all travelers in the same origin-destination (O-D) pair are equal, and no traveler can lower his/her cost by unilaterally switching to a different route. 

\item Wardrop's second principle, known as the {\it system optimal} principle, assumes that travelers behave cooperatively in making their travel decisions such that the total travel cost on the entire network is minimized. In this case, the travel costs experienced by travelers in the same O-D pair are not necessarily identical.

\end{itemize}

Since the seminal work of \cite{MNa, MNb}, the DTA literature has been focusing on the dynamic extension of Wardrop's principles, which gives rise to the notions of {\it dynamic user equilibrium} (DUE) and {\it dynamic system optimal} (DSO) models. The DUE model stipulates that experienced travel cost, including travel time and early/late arrival penalties, is identical for those route and departure time choices selected by travelers between a given O-D pair. The DSO model seeks  system-wide minimization of travel costs incurred by all the travelers subject to the constraints of fixed travel demand and network flow dynamics. 

For a comprehensive review of DTA models, the reader is referred to \cite{BLR2001, PZ2001, SL2005, SL2006, J2007} and \cite{BRBBWS2017}. A discussion of DTA from the perspective of intelligent transportation system can be found in \cite{RB1996}. \cite{CBMPBWH2011} present a primer on simulation-based DTA modeling. \cite{GHP2016} focus on the traffic flow modeling aspect of DTA, namely the hydrodynamic models for vehicular traffic and their network extensions. \cite{WSHF2018} review relevant DTA literature concerning environmental sustainability.

Dynamic user equilibrium (DUE), which is one type of DTA, remains a modern perspective on traffic network modeling that enjoys wide scholarly support. It is conventionally studied as open-loop, non-atomic Nash-like games \citep{FBSTW1993}. The notion of open loop refers to the assumption that the travelers' route choices do not change with time or in response to dynamic network conditions after they leave the origin. The non-atomic nature refers to the prevailing technique of flow-based modeling, instead of treating the traffic as individual vehicles; this is in contrast to agent-based modeling \citep{BCNR2004, CBN2003, SHOA2017}. DUE captures two aspects of travel behavior quite well: departure time choice and route choice. Within the DUE model, travel cost for the same trip purpose is identical for all utilized path and departure time choices. The relevant notion of travel cost is a weighted sum of travel time and arrival penalty.

In the last two decades there have been many efforts to develop a theoretically sound formulation of DUE that is also a canonical form acceptable to scholars and practitioners alike. Analytical DUE models tend to be of two varieties: 
\begin{itemize}
\item[(1)] {\it route-choice} (RC) DUE \citep{BB2003, CH1998, LS2002, LHGS2013, RB1996b, TW2000, VD2004, ZM2000}; and 

\item[(2)] {\it simultaneous route-and-departure-time} choice (SRDT) DUE \citep{FBSTW1993, FBST2001, FHNMY2013, FKKR2011, HFY2013b, HFSL2015, HSF2015, HL2002, NZ2010, SL2004, UHD2012, W2002}. 
\end{itemize}

There are two essential components within the RC or SRDT notions of DUE: (1) the mathematical expression of the equilibrium condition; and (2) a network performance model, which is often referred to as {\it dynamic network loading} (DNL). There are multiple means of expressing the Nash-like notion of dynamic equilibrium, including:
\begin{enumerate}
 \setlength\itemsep{0 em}
\item variational inequalities \citep{FBSTW1993, FHNMY2013, HFY2013b, HFSL2015, HSF2015, SW1994, SW1995};

\item nonlinear complementarity problems \citep{HUD2011, PHRU2011, UHD2012, WTC2002};

\item differential variational inequalities \citep{FM2014, FM2006, HFSL2015};

\item differential complementarity system \citep{BPLM2012}; 

\item fixed-point problem in a Hilbert space \citep{FKKR2011, HSF2015}; and

\item stationary point of an evolutionary dynamic \citep{M2006, SW1995}.
\end{enumerate}
On the other hand, the DNL sub-problem captures the relationship between dynamic traffic flows and travel delays, by articulating link dynamics, junction interactions, flow propagation, link delay, and path delay. It has been the focus of a significant number of studies in traffic modeling \citep{RB1996b, XWFMZ1999, TW2000, LS2002, WTC2002, SL2004, YLI2005, PR2006, NZ2010, BPLM2011, UHD2012, HFY2013a}. The DNL model aims at describing and predicting the spatial and temporal evolution of traffic flows on a network that is consistent with established route and departure time choices of travelers. This is done by introducing appropriate dynamics to flow propagation, flow conservation, link delay, and path delay on a network level. Any DNL must be consistent with the established path flows, and is usually performed under the {\it first-in-first-out} (FIFO) rule \citep{SL2006}.

In general, DNL models have the following components:
\begin{enumerate}
 \setlength\itemsep{0 em}
\item some form of link and/or path dynamics;
\item an analytical relationship between flow/speed/density and link traversal time
\item flow propagation constraints; 
\item a model of junction dynamics and delays;
\item a model of path traversal time; and 
\item appropriate initial conditions.
\end{enumerate}

DNL gives rise to the path delay operator, which is analogous to the pay-off function in classical Nash games, and plays a pivotal role in DTA and DUE problems. The properties of the delay operator are critical to the existence and computation of DUE models. However, it is widely recognized that the DNL model or the delay operator is not available in closed form; instead, they have to be numerically evaluated via a computational procedure. As a result, the mathematical properties of the delay operator remain largely unknown. This has significantly impacted the computation of DUE problems due to the lack of provable convergence theories, which all require certain form of generalized monotonicity of the delay operators. Table \ref{tabconv} shows some relevant computational algorithms for DUE and their convergence conditions with respect to the continuity and monotonicity of the delay operators. The reader is referred to \cite{HFSL2015} for definitions of different types of generalized monotonicity.

\begin{table}[h!]
\caption{Computational algorithms for DUE. The algorithms are arranged in an increasing order of generality of the monotonicity.}
\centering
\begin{tabular}{llll}
\hline
                   & {DUE Type}  &  {Computational} & {Convergence} 
\\\hline

  \multirow{2}{*}{\cite{FKKR2011}}  &   \multirow{2}{*}{SRDT DUE} & \multirow{1}{*}{Fixed-point}  & \multirow{1}{*}{Lipschitz cont.}  
 \\
 \multirow{2}{*}{}  &   \multirow{1}{*}{}  &   \multirow{1}{*}{algorithm} &   \multirow{1}{*}{strongly monotone}

 \\\hline

 \multirow{2}{*}{\cite{LS2002}}  &   \multirow{2}{*}{{RC DUE}} & \multirow{1}{*}{{Alternating direction}}  & \multirow{2}{*}{Co-coercive}  
 \\
 \multirow{2}{*}{}  &   \multirow{2}{*}{}  &   \multirow{1}{*}{{algorithm}} &   \multirow{2}{*}{{}}

 \\\hline

 \multirow{2}{*}{\cite{SL2004}}  &   \multirow{1}{*}{{SRDT DUE}} & \multirow{1}{*}{descent algorithm}  & \multirow{2}{*}{Co-coercive}  
 \\
 \multirow{2}{*}{}  &   \multirow{1}{*}{elastic demand}  &   \multirow{1}{*}{(projection)} &   \multirow{2}{*}{}

 \\\hline

  \multirow{2}{*}{\cite{SL2006}}  &   \multirow{1}{*}{RC DUE} & \multirow{1}{*}{Route-swapping}  & \multirow{1}{*}{Continuous}  
 \\
 \multirow{2}{*}{}  &   \multirow{1}{*}{bounded rationality}  &   \multirow{1}{*}{algorithm} &   \multirow{1}{*}{monotone}

 \\\hline

  \multirow{2}{*}{\cite{HL2002}}  &   \multirow{2}{*}{SRDT DUE} & \multirow{1}{*}{Route-swapping}  & \multirow{1}{*}{Continuous}  
 \\
 \multirow{2}{*}{}  &   \multirow{1}{*}{}  &   \multirow{1}{*}{algorithm} &   \multirow{1}{*}{monotone}

 \\\hline

  \multirow{2}{*}{\cite{THG2012}}  &   \multirow{2}{*}{SRDT DUE} & \multirow{1}{*}{Route-swapping}  & \multirow{1}{*}{Continuous}  
 \\
 \multirow{2}{*}{}  &   \multirow{1}{*}{}  &   \multirow{1}{*}{algorithm} &   \multirow{1}{*}{monotone}

 \\\hline

\multirow{2}{*}{ \cite{LHGS2013}} & \multirow{2}{*}{{RC DUE}}  & \multirow{1}{*}{{Extragradient/}}  &  \multirow{1}{*}{{Lipschitz cont.}} 
\\
 \multirow{2}{*}{}  & \multirow{2}{*}{} & \multirow{1}{*}{double projection}  &\multirow{1}{*}{{pseudo monotone}}

 \\\hline

  \multirow{2}{*}{\cite{HSF2015}}  &   \multirow{1}{*}{SRDT DUE} & \multirow{1}{*}{Self-adaptive}  & \multirow{1}{*}{Continuous}  
 \\
 \multirow{2}{*}{}  &   \multirow{1}{*}{bounded rationality}  &   \multirow{1}{*}{projection} &   \multirow{1}{*}{D-property}

 \\\hline

  \multirow{2}{*}{\cite{HFSL2015}}  &   \multirow{1}{*}{SRDT DUE} & \multirow{1}{*}{Proximal point}  & \multirow{2}{*}{Dual solvable}  
 \\
 \multirow{2}{*}{}  &   \multirow{1}{*}{elastic demand}  &   \multirow{1}{*}{method} &   \multirow{2}{*}{}

 \\\hline
\end{tabular}
\label{tabconv}
\end{table}

This paper documents theoretical and numerical advances in synthesizing traffic flow theory and traffic assignment models, by presenting a computational theory of DUE, which includes algorithms and software implementation. While there have been numerous studies on the modeling and computation of DUEs, including those reviewed in this paper, little agreement exists regarding an appropriate mathematical formulation of DUE or DNL models, as well as the extent to which certain models should/can be applied. This is partially due to the lack of open-source solvers and a set of benchmarking test problems for DUE models. In addition, large-scale computational examples of DUE were rarely reported and, when they were, little detail was provided that allows the results to be \emph{validated}, \emph{reproduced} and \emph{compared}. 

This paper aims to bridge the aforementioned gap by presenting a computable theory for the simultaneous route-and-departure-time (SRDT) DUE model along with open-source software packages. In particular, the DNL model is based on the Lighthill-Whitham-Richards fluid dynamic model \citep{LW1955, R1956}, and is formulated as a system of differential algebraic equations (DAEs) by invoking the variational theory for kinematic wave models. This technique allows the DAE system to be solved with straightforward time-stepping without the need to solve any partial differential equations. Moreover, this paper presents the fixed-point algorithm for solving the DUE problem, which was derived from the {\it differential variational inequality} formalism \citep{FH2018}. Both the DNL procedure and the fixed-point algorithm are implemented in MATLAB, and we present the computational results on several test networks, including the Chicago Sketch network with 250,000 paths. To our knowledge, this is the largest instance of SRDT DUE computation reported in the literature to date.

In addition, we make openly available the MATLAB package, which can be used to solve DUEs on user-defined networks. The package is documented in this paper, and the programs are attached to this publication. More updates and examples are available at \textcolor{blue}{https://github.com/DrKeHan/DTA}. It is our intention that the open-source package will not only help DTA modelers with benchmarking a wide range of algorithms and solutions, but also offer researchers a platform to further develop their own models and applications.

The rest of the paper is organized as follows. Section \ref{secFormulation} introduces some key notions and mathematical formulations of DUE. Section \ref{secDNL} details the dynamic network loading procedure and the DAE system formulation. The fixed-point algorithm for computing DUE is presented in Section \ref{secFPalgo}. The computational results obtained from the proposed DUE solver are presented in Section \ref{secResults}. Section \ref{secConclusion} offers some concluding remarks. Finally, the MATLAB software package is documented in the Appendix.

\section{Formulations of dynamic user equilibrium}\label{secFormulation}

 We introduce a few notations and terminologies for the ease of presentation below.

\begin{itemize}
 \setlength\itemsep{0 em}
 \setlength{\itemindent}{-.35 in}
\item[]{\makebox[1.6cm]{$\mathcal{P}$\hfill} set of paths in the network}

\item[]{\makebox[1.6cm]{$\mathcal{W}$\hfill} set of O-D pairs in the network}

\item[]{\makebox[1.6cm]{$Q_{ij}$\hfill}  fixed O-D demand between $(i,\,j)\in\mathcal{W}$}

\item[]{\makebox[1.6cm]{$\mathcal{P}_{ij}$\hfill}  subset of paths that connect O-D pair $(i,\,j)$}

\item[]{\makebox[1.6cm]{$t$\hfill} continuous time parameter in a fixed time horizon $[t_0,\,t_f]$}

\item[]{\makebox[1.6cm]{$h_p(t)$\hfill} departure rate along path $p$ at time $t$}

\item[]{\makebox[1.6cm]{$h(t)$\hfill} complete vector of departure rates $h(t)=\big(h_p(t):\, p\in\mathcal{P}\big)$}

\item[]{\makebox[1.6cm]{$\Psi_p(t,\,h)$\hfill} travel cost along path $p$ with departure time $t$, under departure profile $h$}

\item[]{\makebox[1.6cm]{$v_{ij}(h)$\hfill}  minimum travel cost between O-D pair $(i,\,j)$ for all paths and departure times}

\end{itemize}

\noindent We stipulate that the path departure rates are square integrable: 
$$
h_p(\cdot)\in L_+^2[t_0,\,t_f],\qquad h(\cdot)\in\big(L_+^2[t_0,\,t_f]\big)^{|\mathcal{P}|}
$$
\noindent We define the {\it effective delay operator} $\Psi$ as follows:
\begin{equation}
\begin{array}{c}
\Psi:~\big(L_+^2[t_0,\,t_f]\big)^{|\mathcal{P}|}~\to~\big(L_{+}^2[t_0,\,t_f]\big)^{|\mathcal{P}|}
\\
\\
h(\cdot)=\big\{ h_p(\cdot),\,p\in\mathcal{P}\big\}~\mapsto~\Psi(h)=\big\{\Psi_p(\cdot,\,h),\,p\in\mathcal{P} \big\}
\end{array}
\end{equation}

\noindent Here, the term `effective delay' \citep{FBSTW1993} is a generalized notion of travel cost that may include not only a linear combination of travel time and arrival penalty, but also other forms of cost such as road pricing. The effective delay operator $\Psi$ is essential to the DUE model as it encapsulates the physics of the traffic network by capturing the dynamics of traffic flows at the link, junction, path, and network levels. We will discuss this operator in greater detail in Section \ref{secDNL}.

The travel demand satisfaction constraint is expressed as
\begin{equation}
\sum_{p\in\mathcal{P}_{ij}}\int_{t_0}^{t_f} h_p(t)dt~=~Q_{ij}\qquad \forall (i,\,j)\in\mathcal{W}
\end{equation}

\noindent Therefore, the set of feasible path departure vector can be expressed as
\begin{equation}\label{Lambdadef}
\Lambda=\left\{h\geq 0:~\sum_{p\in\mathcal{P}_{ij}}\int_{t_0}^{t_f} h_p(t)dt=Q_{ij}\quad \forall (i,\,j)\in\mathcal{W}\right\}\subset \big(L^2[t_0,\,t_f]\big)^{|\mathcal{P}|}
\end{equation}

\noindent The following definition of dynamic user equilibrium is first proposed by \cite{FBSTW1993} in a measure-theoretic context.
\begin{definition}{\bf (SRDT DUE)} 
A vector of departures $h^*\in\Lambda$ is a dynamic user equilibrium with simultaneous route and departure time (SRDT) choice if 
\begin{equation}\label{DUEdef}
h_p^*(t)>0,~p\in\mathcal{P}_{ij}~\Longrightarrow~\Psi_p(t,\,h^*)=v_{ij}(h^*)\quad \hbox{a.e.}~ t\in[t_0,\,t_f]
\end{equation}
\end{definition}
\noindent where `a.e.', standing for `for almost every', is a technical term employed by measure-theoretic arguments to indicate that \eqref{DUEdef} only needs to hold up to a subset of $[t_0,\,t_f]$ of zero measure.

\subsection{Variational inequality formulation of DUE}

Using measure-theoretic arguments, \cite{FBSTW1993} establish that a SRDT DUE is equivalent to the following variational inequality under suitable regularity conditions:
\begin{equation}\label{DUEVI}
\left. 
\begin{array}{c}
\hbox{find}~h^*\in\Lambda ~\hbox{such that}
\\
\displaystyle \sum_{p\in\mathcal{P}}\int_{t_0}^{t_f}\Psi_p(t,\,h^*)\big(h_p(t)-h_p^*(t)\big)dt\geq 0
\\
\forall h\in\Lambda
\end{array}
\right\}
\end{equation}

\noindent The VI \eqref{DUEVI} may be written in a more generic form by invoking the inner product $\left< \cdot \right>$ in the Hilbert space $\big(L^2[t_0,\,t_f]\big)^{|\mathcal{P}|}$:
$$
\left<f,\,g\right>\doteq \sum_{p\in\mathcal{P}} \int_{t_0}^{t_f} f_p(t) g_p(t)dt \qquad \forall f,\,g\in \big( L^2[t_0,\,t_f]\big)^{|\mathcal{P}|}
$$
\noindent This leads to the following VI representation of DUE:
$$
\left< \Psi(h^*),\, h-h^*\right> \geq 0\qquad\forall h\in\Lambda
$$

\subsection{Nonlinear complementarity formulation of DUE}

The variational inequality formulation \eqref{DUEVI} is equivalent to the following nonlinear complementarity problem \citep{PHRU2011, UHD2012}:
\begin{align}\label{DUECP1}
0\leq h^*_p(t)  & \perp \Psi_p(t,\,h^*)- v_{ij}(h^*)\geq0\qquad \forall p\in\mathcal{P}_{ij},\,(i,\,j)\in\mathcal{W},\,\hbox{a.e.}~t\in[t_0,\,t_f]
\\
\label{DUECP2}
0\leq v_{ij}(h^*) & \perp \displaystyle  \sum_{p\in\mathcal{P}_{ij}} \int_{t_0}^{t_f} h^*_p(t)dt  - v_{ij}(h^*)\geq 0 \qquad \forall (i,\,j)\in\mathcal{W}
\end{align}
\noindent where $a\perp b$ means that $a\cdot b= 0$.

\subsection{Differential variational inequality formulation of DUE}

It is noted in \cite{FKKR2011} that the VI formulation \eqref{DUEVI} of DUE is equivalent to a differential variational inequality (DVI). This is most easily seen by noting that the demand satisfaction constraints may be re-stated as
\begin{equation}
\left.
\begin{array}{l}
\displaystyle {d \over dt}y_{ij}(t)=\sum_{p\in\mathcal{P}_{ij}} h_p(t)\qquad \forall (i,\,j)\in\mathcal{P}_{ij}
\\
y_{ij}(t_0)=0\qquad\forall (i,\,j)\in\mathcal{P}_{ij}
\\
y_{ij}(t_f)=Q_{ij}\qquad \forall (i,\,j)\in\mathcal{P}_{ij}
\end{array}
\right\},
\end{equation}
\noindent which is recognized as a two-point boundary value problem. The DUE may be consequently expressed as a DVI: 
\begin{equation}\label{DUEDVI}
\left.
\begin{array}{c}
\hbox{find}~h^*\in\Lambda_0 ~\hbox{such that}
\\
\displaystyle \sum_{p\in\mathcal{P}}\int_{t_0}^{t_f} \Psi_p(t,\,h^*)\big(h_p(t)-h_p^*(t)\big)dt\geq 0
\\
\forall h\in\Lambda_0
\end{array}
\right\}
\end{equation}
\noindent where 
\begin{equation}\label{Lambda0def}
\Lambda_0=\left\{ h\geq 0:~{d\over dt} y_{ij}(t)=\sum_{p\in\mathcal{P}_{ij}}h_p(t),~y_{ij}(t_0)=0,~y_{ij}(t_f)=Q_{ij}\quad \forall (i,\,j)\in\mathcal{W} \right\}
\end{equation}

The equivalence of DUE and DVI \eqref{DUEDVI} can be shown with elementary optimal control theory applied to a linear-quadratic problem as demonstrated in \cite{FKKR2011}. The DVI formulation is significant because it allows the still emerging theory of differential variational inequalities \citep{DODG, PS2008} to be employed for the analysis and computation of solutions of the DUE problem when simultaneous departure time and route choices are within the purview of users \citep{FM2014, FM2006, HFSL2015}. However, the emerging literature on abstract differential variational inequalities has not been well exploited for either modeling or computing simultaneous route-and-departure-time choice equilibria; this gap was recently bridged by \cite{FH2018}.

\subsection{Fixed-point formulation of DUE}

Experience with variational inequalities suggest that there exists a fixed-point re-statement of the DUE problem in a proper functional space. The fixed-point formulation of continuous-time DUE is first articulated in \cite{FKKR2011} using optimal control theory. We define $P_{\Lambda_0}[\cdot]$ to be the minimum-norm projection operator in the space $\big(L^2[t_0,\,t_f]\big)^{|\mathcal{P}|}$. Then the following fixed-point problem is equivalent to the DUE problem:
\begin{equation}\label{FPformulation}
h^*=P_{\Lambda}\big[ h^* - \alpha\Psi(h^*) \big] \quad \text{or}\quad h^*=P_{\Lambda_0}\big[ h^* - \alpha\Psi(h^*) \big] 
\end{equation}
where $\alpha>0$ is a fixed constant. The equivalence result may be proven by applying the minimum principle and associated optimality conditions to the minimum-norm problem intrinsic to the projection operator $P_{\Lambda_0}[\cdot]$. See \cite{DODG} for the details suppressed here.

\section{Dynamic network loading}\label{secDNL}

An essential component of the DUE formulation is the effective delay operator $\Psi$, which is constructed using the {\it dynamic network loading} (DNL) procedure. This section details one type of DNL models that is based on the fluid dynamic approximation of traffic flow on networks, known as the Lighthill-Whitham-Richards (LWR) model \citep{LW1955, R1956}. The model, including its various discrete forms \citep{N1993a, Daganzo1994, Daganzo1995, YLI2005}, is widely used in the DTA literature. The rest of this section presents a complete DNL procedure based on the LWR model and its variational representation, which lead to a differential algebraic equation (DAE) system.

\subsection{The Lighthill-Whitham-Richards link model} 
The LWR model is capable of describing the physics of kinematic waves (e.g. shock waves, rarefaction waves), and allows network extensions that capture the formation and propagation of vehicle queues as well as vehicle spillback. 

The LWR model describes the spatial and temporal evolution of vehicle density $\rho(t,\,x)$ on a road link using the following partial differential equation:
\begin{equation}\label{LWRPDE}
\partial_t \rho(t,\,x)+\partial_x f\big(\rho(t,\,x)\big)=0\qquad x\in[a,\,b],~ t\in[t_0,\,t_f],
\end{equation}
\noindent where the link of interest is represented as a spatial interval $[a,\,b]$. The fundamental diagram $f(\cdot)$ is continuous, concave, and satisfies $f(\rho)=f(\rho^{\text{jam}})=0$ where $\rho^{\text{jam}}$ denotes the jam density.  Furthermore, there exists a unique critical density value $\rho^{\text{c}}$ where $f(\cdot)$ attains its maximum $f(\rho^{\text{c}})=C$ were $C$ denotes the flow capacity of the link.

A few widely adopted forms of $f(\cdot)$ include the Greenshields \citep{G1935}, the trapezoidal \citep{Daganzo1994, Daganzo1995}, and the triangular \citep{N1993a, N1993b, N1993c} fundamental diagrams. In the remainder of the paper (and also in the MATLAB package), we focus on the following triangular fundamental diagram:
\begin{equation}\label{trifd}
f(\rho)=
\begin{cases}
v\rho\qquad & \rho\in[0,\,\rho^{\text{c}}]
\\
-w (\rho-\rho^{\text{jam}})\qquad & \rho\in(\rho^{\text{c}},\,\rho^{\text{jam}}]
\end{cases}
\end{equation}
\noindent where $v>0$ and $w>0$ denote the forward and backward kinematic wave speeds, respectively.

While \eqref{LWRPDE}  captures the within-link dynamics, the inter-link propagation of congestion requires a careful treatment of junction dynamics, which is underpinned by the notions of link demand and supply.

\subsection{Link demand and supply}\label{subsecLDS}

We consider a road junction with $m$ incoming links and $n$ outgoing links. The dynamic on each of the $m+n$ links is governed by the LWR model \eqref{LWRPDE}; yet these $m+n$ equations are coupled via their relevant boundary conditions. In particular, the following flow conservation constraint must hold:
\begin{equation}\label{fccjunc}
\sum_{i=1}^m f_i\big(\rho_i(t,\,b_i)\big)=\sum_{j=1}^n f_j\big(\rho_j(t,\, a_j)\big)\qquad\forall t\in[t_0,\,t_f],
\end{equation}
\noindent where, without causing any confusion, we always use subscript $i$ or $j$ to indicate the association with link $i$ or $j$. \eqref{fccjunc} simply means that the total flow through the junction is conserved. However, this condition alone does not guarantee a unique flow profile at these $m+n$ links, and additional conditions need to be imposed \citep{GHP2016}. To this end, we define the link demand and supply \citep{LK1999}, where the demand (supply) is viewed as a function of the density near the exit (entrance) of the link:
\begin{align}
& D\big(\rho(t,\,b-)\big)=
\begin{cases}
f\big(\rho(t,\,b-)\big) \qquad & \rho(t,\,b-)<\rho^c
\\
C \qquad & \rho(t,\,b-)\geq \rho^c
\end{cases}
\\
& S\big(\rho(t,\,a+)\big)=
\begin{cases}
C \qquad & \rho(t,\,a+)<\rho^c
\\
f\big(\rho(t,\,a+)\big)\qquad & \rho(t,\,a+)\geq \rho^c
\end{cases}
\end{align}
\noindent Intuitively, the demand (supply) indicates the maximum flow that can exit (enter) the link. That is,
\begin{equation}\label{dsconst}
f_i\big(\rho_i(t,\,b_i)\big)\leq D_i\big(\rho_i(t,\,b_i-)\big),\quad f_j\big(\rho_j(t,\,a_j)\big)\leq S_j\big(\rho_j(t,\,a_j+)\big)
\end{equation}
\noindent for $i\in\{1,\ldots, m\}$, $j\in\{1,\ldots, n\}$. Similar to \eqref{fccjunc}, \eqref{dsconst} ensures the physical feasibility of the flows through the junction. Nevertheless, additional conditions are needed to isolate a unique flow profile at this junction; these conditions are often derived based on driving behavior or traffic management measures, such as flow distribution \citep{Daganzo1995}, right of way \citep{Daganzo1995, CGP2005}, and traffic signal control \citep{HGPFY2014, HG2015}. A review of these different models is provided in \cite{GHP2016}.

\subsection{The variational representation of link dynamics}
The variational solution representation of Hamilton-Jacobi equations has been widely applied to the LWR-based traffic modeling \citep{N1993a, N1993b, N1993c, Daganzo2005, Daganzo2006, CB2010, LL2013, CL2014, HYJF2017}. We consider the Moskowitz function \citep{Moskowitz1965}, $N(t,\,x)$, which measures the cumulative number of vehicles that have passed location $x$ along a link by time $t$. The following identities hold:
$$
\partial_t N(t,\,x)=f\big(\rho(t,\,x)\big),\quad \partial_x N(t,\,x)=-\rho(t,\,x).
$$
\noindent It is easy to show that $N(t,\,x)$ satisfies the following Hamilton-Jacobi equation:
\begin{equation}
\partial_t N(t,\,x)-f\left(-\partial_x N(t,\,x)\right)=0\qquad x\in[a,\,b],~t\in[t_0,\,t_f].
\end{equation}
\noindent Next, we denote by $f^{\text{in}}(t)$ and $f^{\text{out}}(t)$ the link inflow and outflow, respectively. The cumulative link entering and exiting vehicle counts are defined as
$$
N^{\text{up}}(t)= \int_{t_0}^t f^{\text{in}}(s)ds,\qquad N^{\text{dn}}(t)=\int_{t_0}^t f^{\text{out}}(s)ds,
$$ 

\noindent where the superscripts `up' and `dn' represent the {\it upstream} and {\it downstream} boundaries of the link, respectively. \cite{HPS2016} derive explicit formulae for the link demand and supply based on a variational formulation known as the Lax-Hopf formula \citep{ABS2008, CB2010}, as follows:
\begin{align}\label{DSvteqn1}
& D(t)=
\begin{cases}
f^{\text{in}}\left(t- {L\over v}\right) \qquad & \text{if}~ N^{\text{up}}\left(t-{L\over v}\right)=N^{\text{dn}}(t)
\\
C \qquad & \text{if}~ N^{\text{up}}\left(t-{L\over v}\right)>N^{\text{dn}}(t)
\end{cases}
\\
\label{DSvteqn2}
& S(t)=
\begin{cases}
f^{\text{out}}\left(t- {L\over w}\right) \qquad & \text{if}~ N^{\text{up}}(t)=N^{\text{dn}}\left(t-{L\over w}\right)+\rho^{\text{jam}}L
\\
C \qquad & \text{if}~ N^{\text{up}}(t)<N^{\text{dn}}\left(t-{L\over w}\right)+\rho^{\text{jam}}L
\end{cases}
\end{align}
\noindent where $L=b-a$ denotes the link length. Note that \eqref{DSvteqn1} and \eqref{DSvteqn2} express the link demand and supply, which are inputs of the junction model, in terms of $N^{\text{up}}(\cdot)$ and $N^{\text{dn}}(\cdot)$, or $f^{\text{in}}(\cdot)$ and $f^{\text{out}}(\cdot)$. This means that one no longer needs to compute the dynamics within the link, but to focus instead on the flows or cumulative counts at the two boundaries of the  link. Such an observation tremendously simplifies the link dynamics, and gives rise to the link-based formulation \citep{Jin2015, HPS2016} or, in its discrete form, the link transmission model \citep{YLI2005}.

\subsection{Junction dynamics that incorporate route information}

Essential to the network extension of the LWR model is the junction model. Unlike many existing junction models such as those reviewed in Section \ref{subsecLDS}, in a path-based DNL procedure one must incorporate established routing information into the junction model. Such information is manifested in an endogenous {\it flow distribution matrix}, which specifies the portion of exit flow from a certain incoming link that advances into a given outgoing link. This can be done by explicitly tracking the route composition in every unit of flow along the link.

We begin by defining the link entry time function $\tau(t)$ where $t$ denotes the exit time. Such a function can be obtained by evaluating the horizontal difference between the cumulative curves $N^{\text{up}}(\cdot)$ and $N^{\text{dn}}(\cdot)$ \citep{FHNMY2013}. Next, for link $i$ and path $p$ such that $i\in p$, we define $\mu_i^p(t,\,x)$ to be the percentage of flow on link $i$ that belongs to path $p$. The first-in-first-out principle yields the following identity:
\begin{equation}
\mu_i^p(t,\,b_i)=\mu_i^p\big(\tau_i(t),\,a_i\big)
\end{equation}

Now we consider a junction $J$ with incoming links labeled as $i\in\{1,\ldots, m\}$ and outgoing links labeled as $j\in\{1,\ldots, n\}$. The distribution matrix $A^J(t)$ can be expressed as
$$
A^J(t)=\big\{\alpha_{ij}(t) \big\},\qquad \alpha_{ij}(t)=\sum_{p\ni i,j} \mu_i^p\big(\tau_i(t),\,a_i\big)
$$
\noindent There exist a number of choices for the junction models; they all need to satisfy the flow conservation constraint \eqref{fccjunc}, the demand-supply constraints \eqref{dsconst}, and depend on the flow distribution matrix $A^J(t)$. Any such model can be conceptually expressed as:
\begin{equation}\label{RSconceptual}
\left(\big[f_i^{\text{out}}(t)\big]_{i=1,\ldots, m},\, \big[f_j^{\text{in}}(t)\big]_{j=1,\ldots, n} \right)=\Theta\left(\big[D_i(t)\big]_{i=1,\ldots,m},\, \big[S_j(t)\big]_{j=1,\ldots,n}; \, A^J(t)\right)
\end{equation}
\noindent where $\Theta$ denotes the junction model where $D_i(t)$, $S_j(t)$ and $A^J(t)$ are treated as its input arguments. The output of the model, shown as the right hand side of \eqref{RSconceptual}, include the outflows (inflows) of the incoming (outgoing) links.

\subsection{Dynamics at the origin nodes}

A model at the origin (source) nodes is needed since the path flows $h_p(\cdot)$, defined by \eqref{Lambdadef}, are not bounded from above. In this case, a queuing model is needed at the origin node in case the relevant departure rate exceed the flow capacity of the first link. 

We employ a simple point-queue type dynamic \citep{V1969} for the origin node $o$. Denote by $q_o(t)$ the volume of the point queue, and let link $j$ be the link connected to the origin node. We have that
\begin{equation}\label{sourcpqm}
{d\over dt} q_o(t)=\sum_{p\in \mathcal{P}^o}h_p(t) -  \min\big\{D_o(t),\, S_j(t)\big\}
\end{equation}
\noindent where $\mathcal{P}^o$ denotes the set of paths originating from $o$. The first term on the right hand side of \eqref{sourcpqm} represents flow into the point queue, while the second term represents flow leaving the queue, where the demand at the origin is defined as
$$
D_o(t)=
\begin{cases}
\mathcal{M} \qquad & q_o(t)>0
\\\\
\displaystyle \sum_{p\in \mathcal{P}^o}h_p(t) \qquad & q_o(t)=0
\end{cases},
$$
\noindent and $\mathcal{M}$ is a sufficiently large number, e.g. larger than the flow capacity of link $j$.

\subsection{Calculating path travel times}

The DNL procedure calculates the path travel times with given path departure rates. The path travel time consists of link travel times plus possible queuing time at the origin. We define the link exit time function $\lambda(t)$ by measuring the horizontal difference between the cumulative entering and exiting counts:
\begin{equation}
N^{\text{up}}(t)=N^{\text{dn}}\big(\lambda(t)\big)
\end{equation}
\noindent For a path expressed as $p=\{1,\,2,\,\ldots,\,K\}$, the path travel time $D_p(t,\,h)$ is calculated as
\begin{equation}
\lambda_s \circ \lambda_1 \circ \lambda_2 \ldots \circ \lambda_{K}(t)
\end{equation}
\noindent where $f\circ g(t)\doteq g\big(f(t)\big)$ denotes the composition of two functions. $\lambda_o(\cdot)$ is the exit time function for the potential queuing at the origin $o$.

\subsection{The differential algebraic equation system formulation of DNL}

Here, as a summary of the individual sections presented so far, we present the complete  differential algebraic equation (DAE) system formulation. We begin with the following list of key notations. 

\begin{itemize}
 \setlength\itemsep{0 em}
 \setlength{\itemindent}{-.35 in}
\item[]{\makebox[1.6cm]{$\mathcal{P}$\hfill}  set of all paths}

\item[]{\makebox[1.6cm]{$\mathcal{S}$\hfill}  set of origins}

\item[]{\makebox[1.6cm]{$\mathcal{P}^o$\hfill}  set of paths originating from $o\in\mathcal{S}$}

\item[]{\makebox[1.6cm]{$\mathcal{I}^J$\hfill}  set of incoming links of a junction $J$}

\item[]{\makebox[1.6cm]{$\mathcal{O}^J$\hfill}  set of outgoing links of a junction $J$}

\item[]{\makebox[1.6cm]{$A^J$\hfill}  flow distribution matrix of junction $J$}

\item[]{\makebox[1.6cm]{$h_p(t)$\hfill}  departure rate along path $p\in\mathcal{P}$}

\item[]{\makebox[1.6cm]{$f^{\text{in}}_i(t)$\hfill}  inflow of link $i$}

\item[]{\makebox[1.6cm]{$f^{\text{out}}_i(t)$\hfill}  outflow of link $i$}

\item[]{\makebox[1.6cm]{$N^{\text{up}}_i(t)$\hfill} cumulative link entering count}

\item[]{\makebox[1.6cm]{$N^{\text{dn}}_i(t)$\hfill} cumulative link exiting count}

\item[]{\makebox[1.6cm]{$\mu_i^p(t,\,x)$\hfill}  percentage of flow on link $i$ that belongs to path $p$}

\item[]{\makebox[1.6cm]{$q_o(t)$\hfill} point queue at the origin node $o\in\mathcal{S}$}

\item[]{\makebox[1.6cm]{$\tau_i(t)$\hfill}  entry time of link $i$ corresponding to  exit time $t$}

\item[]{\makebox[1.6cm]{$\lambda_i(t)$\hfill}  exit time of link $i$ corresponding to  entry time $t$}

\end{itemize}

\noindent The DAE system reads:
\begin{align}
\label{DAEeqn1}
& {d\over dt} q_o(t)=\sum_{p\in \mathcal{P}^o}h_p(t) -  \min\big\{D_o(t),\, S_j(t)\big\},\quad 
D_o(t)=
\begin{cases}
\mathcal{M} \qquad & q_o(t)>0
\\
 \sum_{p\in \mathcal{P}^o}h_p(t) \qquad & q_o(t)=0
\end{cases}
\\
\label{DAEeqn2}
& D_i(t)=
\begin{cases}
f_i^{\text{in}}\left(t- {L_i\over v_i}\right) \quad & \text{if}~ N_i^{\text{up}}\left(t-{L_i\over v_i}\right)=N_i^{\text{dn}}(t)
\\
C_i \quad & \text{if}~ N_i^{\text{up}}\left(t-{L_i\over v_i}\right)>N_i^{\text{dn}}(t)
\end{cases}
\\
\label{DAEeqn3}
& S_j(t)=
\begin{cases}
f_j^{\text{out}}\left(t- {L_j\over w_j}\right) \quad & \text{if}~ N_j^{\text{up}}(t)=N_j^{\text{dn}}\left(t-{L_j\over w_j}\right)+\rho_j^{\text{jam}}L_j
\\
C_j \quad & \text{if}~ N_j^{\text{up}}(t)<N_j^{\text{dn}}\left(t-{L_j\over w_j}\right)+\rho_j^{\text{jam}}L_j
\end{cases}
\\
\label{DAEeqn4}
& N_i^{\text{dn}}(t)=N_i^{\text{up}}\big(\tau_i(t)\big),\qquad N_i^{\text{up}}(t)=N_i^{\text{dn}}\big(\lambda_i(t)\big)
\\
\label{DAEeqn5}
& \mu_j^p(t,\,a_j)={f_i^{\text{out}}(t)\mu_i^p\big(\tau_i(t),\,a_i\big)\over f_j^{\text{in}}(t)} \quad \forall p ~\text{s.t.}~\{i,\,j\}\subset p
\\
\label{DAEeqn6}
& A^J(t)=\big\{\alpha_{ij}(t)\big\},\quad \alpha_{ij}(t)=\sum_{p\ni i, j} \mu_i^p\big(\tau_i(t),\,a_i\big)
\\
\label{DAEeqn7}
& \left(\big[f_i^{\text{out}}(t+)\big]_{i=1,\ldots, m},\, \big[f_j^{\text{in}}(t+)\big]_{j=1,\ldots, n} \right)=\Theta\left(\big[D_i(t)\big]_{i=1,\ldots,m},\, \big[S_j(t)\big]_{j=1,\ldots,n}; \, A^J(t)\right)
\\
\label{DAEeqn8}
& {d\over dt} N_i^{\text{up}}(t)=f_i^{\text{in}}(t),\quad {d\over dt} N_i^{\text{dn}}(t)=f_i^{\text{out}}(t)
\\
\label{DAEeqn9}
& D_p(t,\,h)=\lambda_s \circ \lambda_1 \circ \lambda_2 \ldots \circ \lambda_{K}(t)\qquad p=\{1,\,2,\,\ldots,\,K\}
\end{align}

\noindent Eqns \eqref{DAEeqn1}-\eqref{DAEeqn9} form the DAE system for the DNL procedure. Compared to the partial differential algebraic equation (PDAE) system presented in \cite{HPF2016a}, the DAE system does not involve any spatial derivative as one would expected from the LWR-type equations, by virtue of the variational formulation.

\begin{figure}[h!]
\centering
\includegraphics[width=\textwidth]{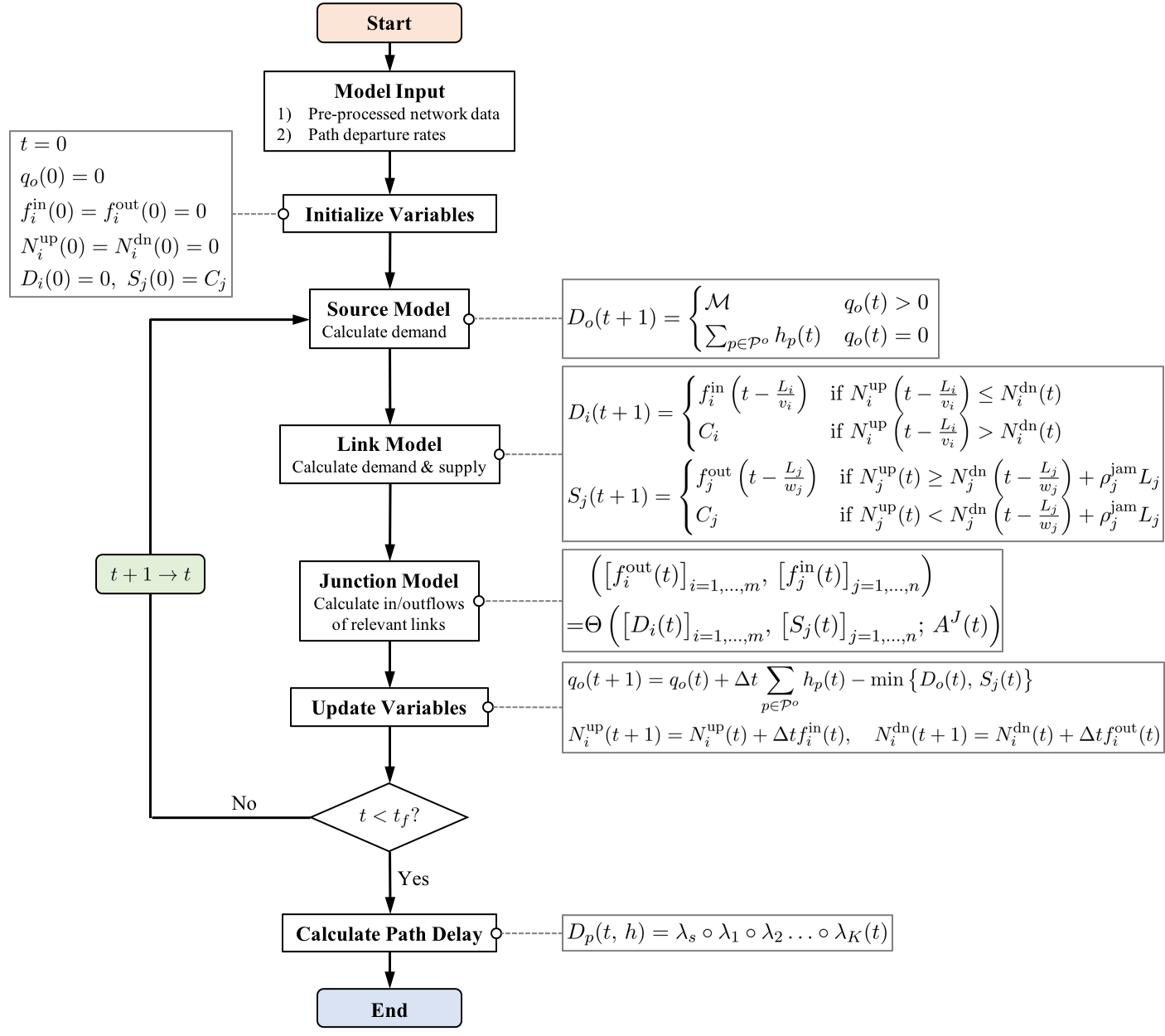}
\caption{Time-stepping logic of the discretized DAE system \eqref{DAEeqn1}-\eqref{DAEeqn9}. Here, $t=0$, $1$, $2$, $\ldots$, represents discrete time steps with step size denoted $\Delta t$.}
\label{figFlow}
\end{figure}

The proposed DAE system may be time-discretized and solved in a forward fashion. Figure \ref{figFlow} explains the time-stepping logic, which is implemented in the Matlab package associated with this publication.

\section{The fixed-point algorithm for computing DUE}\label{secFPalgo}

The computation of DUE is facilitated by the equivalent mathematical formulations such as variational inequality, differential variational inequality, fixed-point problem, and nonlinear complementarity problem. Some of these methods and their convergence conditions are mentioned in Table \ref{tabconv}. In this section, we present the following fixed-point algorithm based on the fixed-point formulation \eqref{FPformulation}:
\begin{equation}\label{FPiterate}
h^{k+1}=P_{\Lambda_0}\left[h^k-\alpha\Psi(h^k)\right]
\end{equation}
\noindent where $\alpha>0$ is a constant, $h^{k+1}$ and $h^{k}$ respectively represent the path departure rate vector at the $k+1$- and $k$-th iteration; $\Psi(h^k)$ denote the effective delays. Recalling the definition of $\Lambda_0$ from \eqref{Lambda0def}, the right hand side of \eqref{FPiterate} amounts to a linear quadratic optimal control problem whose dual variables can be explicitly found following Pontryagin's minimum principle; the reader is referred to \cite{FKKR2011} for details. 

We outline the main steps of the fixed-point algorithm below.

\begin{framed}
\noindent {\bf Fixed-Point Algorithm for Solving DUE}
\begin{description}
\item \textbf{Step 0. Initialization}. Set $k=0$ and select an initial
departure rate vector $h^{0}\in\Lambda$. Fix a suitable constant $\alpha >0$ used in all iterations. 

\item \textbf{Step 1. Dynamic Network Loading}. Carry out a dynamic network loading procedure with departure rate vector $h^k\in\Lambda$ to compute the effective path delays $\Psi_p(t,\,h^k)$ for all $t\in[t_0,\,t_f]$ and $p\in\mathcal{P}$.

\item \textbf{Step 2. Fixed-Point Update.} For every origin-destination pair $(i,\,j)\in \mathcal{W}$, solve the following algebraic equation for the dual variable $v_{ij}$ (where $[x]_{+}\doteq \max \{0,\,x\}$ assures non-negativity): 
\begin{equation}\label{eqnFPalgupdate}
\sum_{p\in \mathcal{P}_{ij}}\int_{t_{0}}^{t_{f}}\left[ h_{p}^{k}(t)-\alpha
\Psi _{p}(t,\,h^{k})+v_{ij}\right] _{+}dt=Q_{ij}
\end{equation}
\noindent For all $t\in[t_0,\,t_f]$ and $p\in\mathcal{P}_{ij}$, compute
\begin{equation*}
h_p^{k+1}(t) =\left[ h_p^{k}(t)-\alpha \Psi_p(t,\,h^{k})+v_{ij}\right]_{+}
\end{equation*}

\item \textbf{Step 3. Stopping Test}. For a predetermined tolerance $\epsilon>0$, if
\begin{equation*}
{\|  h^{k+1}-h^{k} \|^2 \over \| h^k\|^2} \leq \epsilon
\end{equation*}
\noindent stop and declare $h^{k+1}$ a DUE solution. Otherwise set $k=k+1$ and go to {\bf Step 1}.
\end{description}
\end{framed}

\begin{remark} In the fixed-point algorithm, the critical step is to find the dual variable $v_{ij}$ 
in \eqref{eqnFPalgupdate}. Note that this amounts to finding $x$ such that $G(x)=0$, where 
$$
G(x)\doteq \sum_{p\in\mathcal{P}_{ij}}\int_{t_0}^{t_f} \left[h_{p}^{k}(t)-\alpha
\Psi _{p}(t,\,h^{k})+ x \right] _{+}dt-Q_{ij}
$$
\noindent is a continuous function with a single argument $x$. Therefore, the dual variable can be found via standard root-finding algorithms such as Bracketing or Bisection methods. 
\end{remark}

\section{Computational examples of DNL and DUE}\label{secResults}

\begin{figure}[h!]
\centering
\includegraphics[width=.8\textwidth]{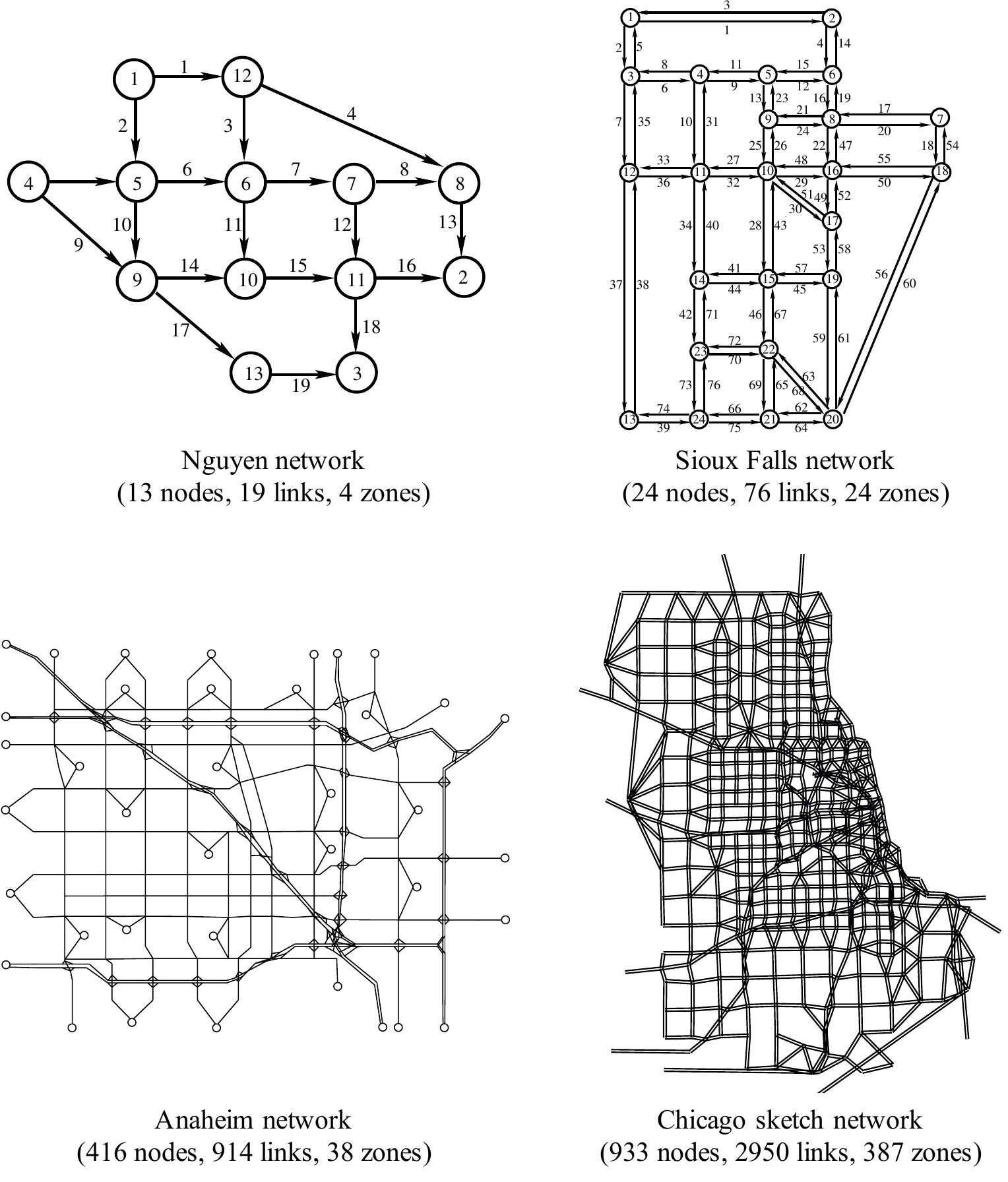}
\caption{The four test networks for DUE algorithms.}
\label{figchapNum:Fournetworks}
\end{figure}

We present computational examples of the simultaneous route-and-departure-time dynamic user equilibria on four networks of varying sizes and shapes, as shown in Table \ref{chapNum:tabfournetworks} and Figure \ref{figchapNum:Fournetworks}. In particular, the Nguyen network was initially studied in \cite{Nguyen1984}, and the last three networks are based on real-world cities in the US, although different levels of network aggregation and simplifications have been applied. Detailed network parameters, including coordinates of nodes and link attributes, are sourced and adapted from \cite{BarGera2018}. Given that the DUE and DNL formulations in this paper are path-based, enumeration of paths was applied to generate path set using the Frank-Wolfe algorithm.

\begin{table}[h!]
\centering
\begin{tabular}{l|l|l|l|l}
   &  Nguyen network  &  Sioux Falls  & Anaheim  & Chicago Sketch
\\\hline
No. of links  &  19 &  76 & 914 &   2,950
\\
No. of nodes  & 13 & 24  &  416  &  933
\\
No. of zones   & 4 & 24  &  38  &  387
\\
No. of O-D pairs  & 4 & 530  &  1,406  & 86,179
\\
No. of paths   & 24 &  6,180  &  30,719  &  250,000
\end{tabular}
\caption{Key attributes of the test networks.}
\label{chapNum:tabfournetworks}
\end{table}

We apply the fixed-point algorithm (Section \ref{secFPalgo}) with the embedded DNL procedure (Section \ref{secDNL}). The fixed-point algorithm is chosen among many other alternatives in the literature, as our extensive experience with DUE computations suggests that this method tends to exhibit satisfactory empirical convergence within limited number of iterations, despite that its theoretical convergence requires strong monotonicity of the delay operator. The DNL sub-model is solved as a DAE system \eqref{DAEeqn1}-\eqref{DAEeqn9}, following the time stepping logic in Figure \ref{figFlow}.

All the computations reported in this section were performed using the MATLAB (R2017b)  package on a standard desktop with Intel i5 processor and 8 GB of RAM.

\subsection{Performance of the fixed-point algorithm}

\begin{figure}[h!]
\centering
\includegraphics[width=\textwidth]{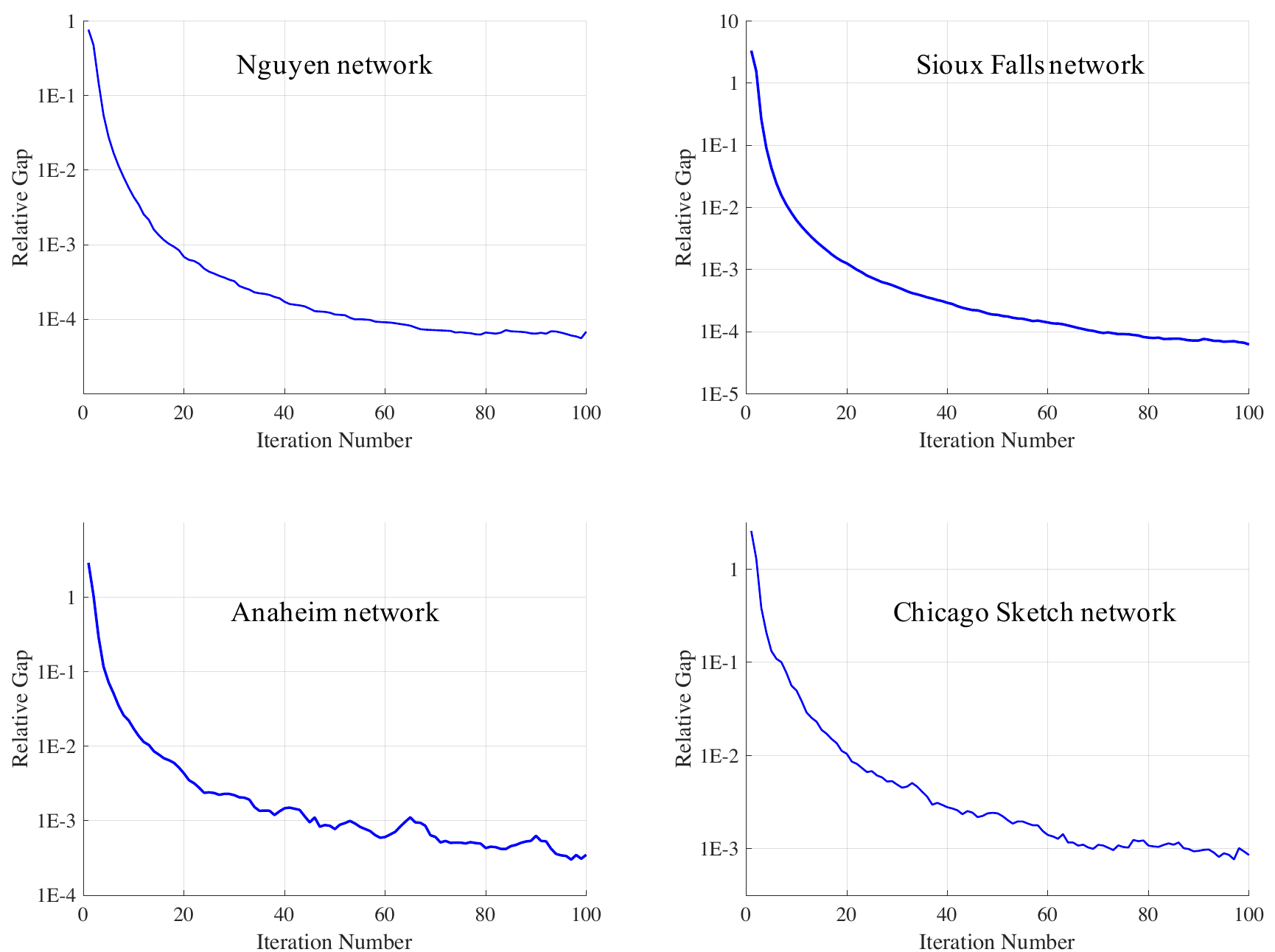}
\caption{The relative gaps (base 10 logarithm) within 100 fixed-point iterations. Here, 1E-x means $10^{-x}$.}
\label{figchapNum:Convergence100}
\end{figure}

The termination criterion for the fixed-point algorithm is set as follows:
\begin{equation}\label{chapNum:fptermination}
{\left\| h^{k+1}-h^{k}\right\|^2 \over \left\|h^{k}\right\|^2} \leq \epsilon
\end{equation}
\noindent where $h^{k}$ denotes the path departure vector in the $k$-th iteration. The threshold $\epsilon$ is set to be $10^{-4}$ for the Nguyen and Sioux Falls networks, and $10^{-3}$ for the Anaheim and Chicago Sketch networks. These different thresholds were chosen to accommodate the varying convergence performances of the algorithm on different networks (see Figure \ref{figchapNum:Convergence100}).

Table \ref{chapNum:fpalgperformance} summarizes the computational performance of the fixed-point algorithms for different networks based on the termination criterion \eqref{chapNum:fptermination}. It is shown that the same termination criterion requires comparable numbers of iterations for different networks, which suggests the scalability of the algorithms.

\begin{table}[h!]
\centering
\begin{tabular}{l|l|l|l|l}
   &  Nguyen network  &  Sioux Falls  & Anaheim  & Chicago Sketch
\\\hline
No. of iterations  &  54 &  70  &  45 &  69 
\\
Computational time  & 5.9 s & 5.7 min  &  23.3 min  &  4.8 hr
\\
Avg. time per DNL   & 0.1 s & 4.3 s  &  25.7 s  &  163.9 s
\\
Avg. time per FP update  & 0.007 s & 0.6 s  &  3.1 s  & 81.2 s
\end{tabular}
\caption{Performance of the fixed-point algorithm on different networks.}
\label{chapNum:fpalgperformance}
\end{table}

\noindent Figure \ref{figchapNum:Convergence100} shows the relative gaps, i.e. left hand side of \eqref{chapNum:fptermination}, for a total of 100 fixed-point iterations on the four networks. It can be seen that for relatively small networks (Nguyen and Sioux Falls), the convergence can be achieved relatively quickly and to a satisfactory degree; the corresponding curves are monotonically decreasing and smooth. For Anaheim and Chicago Sketch networks, the decreasing trend of the gap can sometimes stall and experience fluctuations locally.

\subsection{DUE solutions}

\begin{figure}[h!]
\centering
\includegraphics[width=\textwidth]{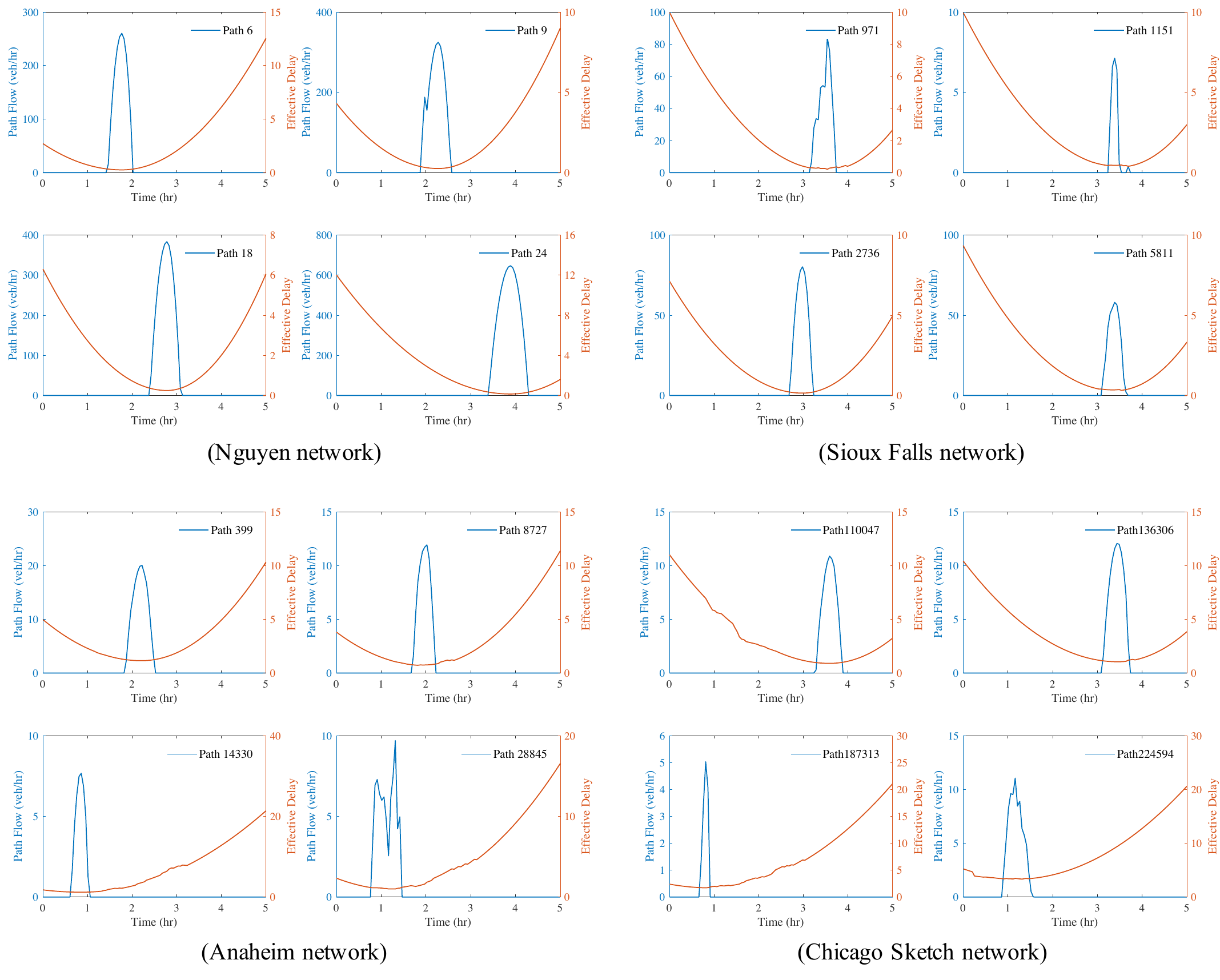}
\caption{Path departure rates and corresponding effective path delays (travel costs) of selected paths in the DUE solutions on the four test networks.}
\label{figDUEsolutions}
\end{figure}

In this section we examine the DUE solutions obtained upon convergence of the fixed-point algorithms. We begin by randomly selecting four paths per network to illustrate the properties of the solutions. Figure \ref{figDUEsolutions} shows the path departure rates as well as the corresponding effective path delays. We observe that the departure rates are non-zero only when the corresponding effective delays are equal and minimum, which conforms to the notion of DUE. Note that the bottoms of the effective delay curves should theoretically be flat, indicating equal travel costs. This is not the case in the figures since we can only obtain approximate DUE solutions in a numerical sense, given the finite number of fixed-point iterations performed to reach those solutions.

\noindent To rigorously assess the quality of the DUE solutions, we define the gap function between each O-D pair $(i,\,j)\in\mathcal{W}$ as
\begin{align}
\text{GAP}_{ij} \doteq & \max\left\{ \Psi_p(t,\,h^*):~ t\in[t_0,\,t_f], p\in\mathcal{P}_{ij} ~\hbox{such that}~h^*_p(t)>0\right\} \nonumber
\\
\label{gapdef}
 & - \min\left\{ \Psi_p(t,\,h^*):~ t\in[t_0,\,t_f], p\in\mathcal{P}_{ij} ~\hbox{such that}~h^*_p(t)>0\right\} 
\end{align}
\noindent Here, $\text{GAP}_{ij}$ represents the range of travel costs experienced by all travelers in the given O-D pair. In an exact DUE solution, the gap should be zero for all O-D pairs.

\begin{figure}[h!]
\centering
\includegraphics[width=\textwidth]{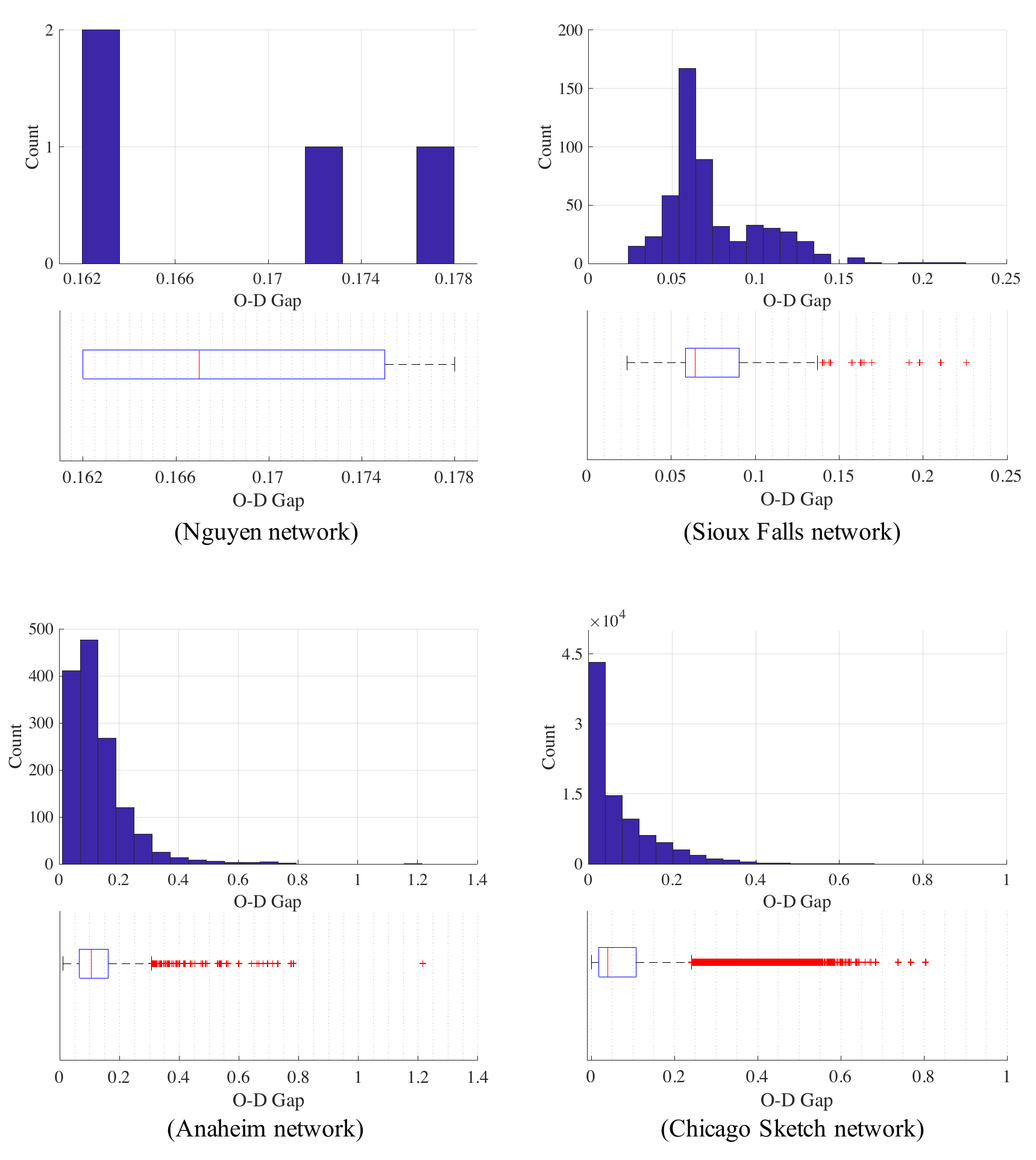}
\caption{Distributions of O-D gaps corresponding to the DUE solutions. The O-D gap is calculated according to \eqref{gapdef}.}
\label{figODgap}
\end{figure}

Figure \ref{figODgap} summarizes all the O-D gaps of the DUE solutions on the four test networks. It can be seen that the majority of the O-D gaps are within 0.2 across all networks. Even for large-scale networks (Anaheim and Chicago Sketch), the 75th percentiles are within 0.15, and the whiskers extend to 0.3. A comparison between Anaheim and Chicago Sketch also reveals that the latter yields better solution quality in terms of O-D gaps, despite the obviously larger size of the problem. This suggests that the solution quality is not necessarily compromised by the size of the network.

\section{Conclusion}\label{secConclusion}

This paper presents a computational theory for dynamic user equilibrium (DUE) on large-scale networks. We begin by presenting a complete and generic dynamic network loading (DNL) model based on the network extension of the LWR model and the variational theory, allowing us to formulate the DNL model as a system of differential algebraic equations (DAEs). The DNL model is capable of capturing the formation, propagation and dissipation of physical queues as well as vehicle spillback; and the DAE system can be discretized and solved in a time-forward fashion. In addition, the fixed-point algorithm for solving DUE problems is reviewed. 

Both the DAE system and the fixed-point algorithm are implemented in MATLAB, and the programs are developed in such a way that they can be applied to solve DUE and DNL problems on any user-defined networks. The MATLAB software package is documented in this paper, which details the structure and flow of the data and individual files.

The MATLAB package is applied to solve DUE problems on several test networks of varying sizes. The largest one is the Chicago Sketch network with 86,179 O-D pairs and 250,000 paths. To the authors' knowledge this is by far the largest instance of SRDT DUE solution reported in the literature, and the codes will be made available along with this publication. Hopefully, our efforts in making these codes and data openly accessible could facilitate the testing and benchmarking of dynamic traffic assignment algorithms, and promote synergies between model development and applications.

\pagebreak

\section*{Appendix: Documentation of the MATLAB software package}

\subsection*{Appendix 1. Dynamic network loading solver}
\begin{figure}[h!]
\centering
\includegraphics[width=.8\textwidth]{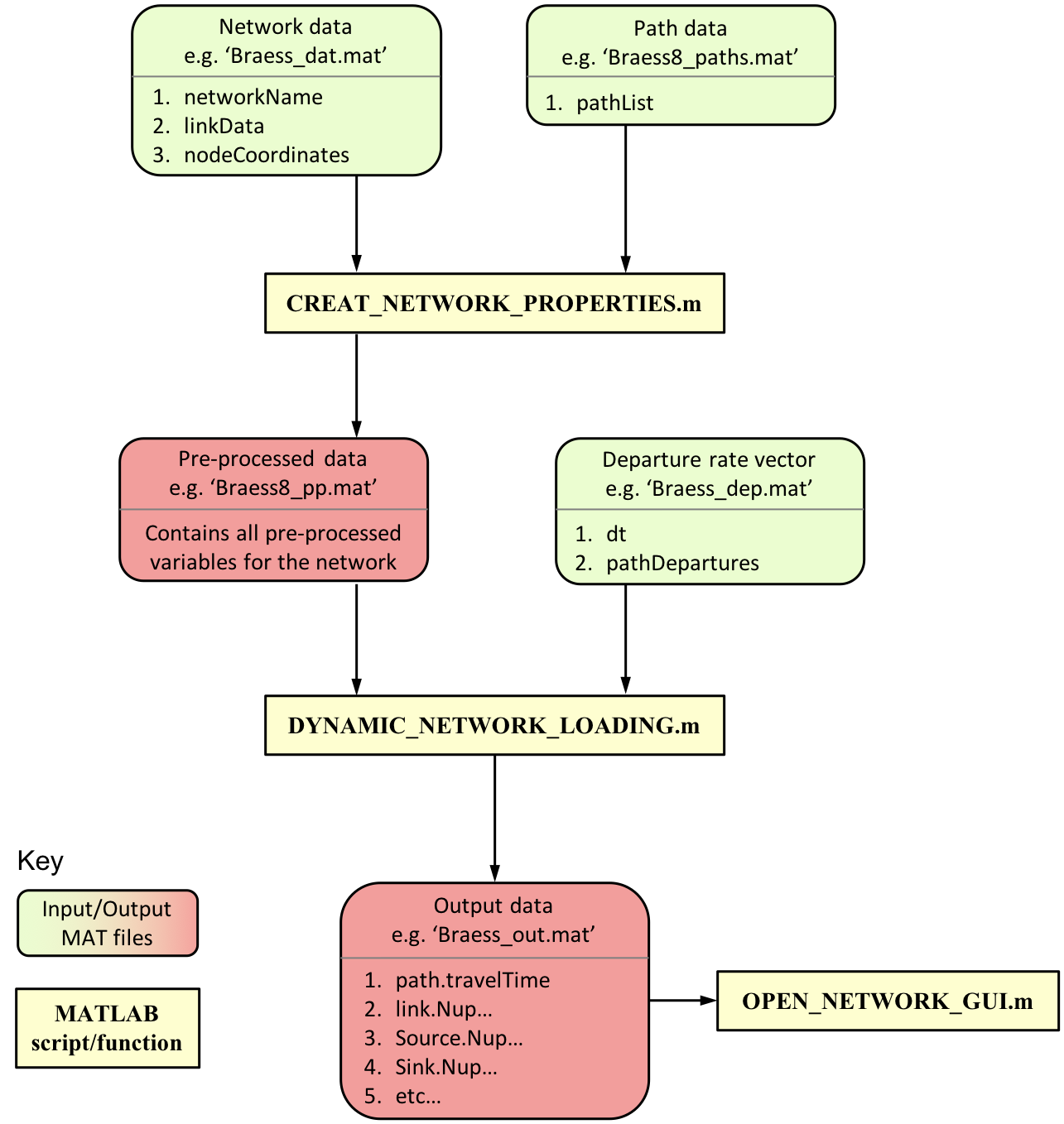}
\caption{Structure and interdependencies of all input and output files within the DNL program. Boxes in green are user-defined data files containing the variables listed.}
\label{figFlow1}
\end{figure}

We document the MATLAB implementation of the dynamic network loading solver, which is based on the discretized DAE system following the procedure in Figure \ref{figFlow}. The individual scripts (.m files) and input data files are illustrated in Figure \ref{figFlow1}.

The DNL solver contains three main script files:\\
\noindent {\bf 1. Pre-processor: `CREATE$\_$NETWORK$\_$PROPERTIES.m'} \\
\noindent The purpose of this script is to prepare all the variables required for the main dynamic network loading script and save them to file or reuse to save computational time. Follow the steps below to run this script.
\begin{itemize}
 \setlength\itemsep{0 em}
\item[1.1.] Set user-defined inputs (at top of file):
\begin{itemize}
 \setlength\itemsep{0 em}
\item[a)] Location of the network data file (e.g. `.../{\bf Braess\_dat.mat}')
\item[b)] Location of path data file (e.g. `.../{\bf Braess\_paths.mat}')
\item[c)] Value of source signal priority\footnote{The DNL code employs junction models that incorporate the continuous signal model proposed by \cite{HGPFY2014} and \cite{HG2015}. For any given junction (node), each of its incoming links, including the source if the node happens to be an origin, will be assigned a priority parameter between 0 and 1, such that the sum of relevant priority parameters equals 1. Once the priority of the source is specified, the priorities of the remaining incoming links are set proportional to the links' capacities. These parameters may be changed dynamically to accommodate different signal control scenarios or strategies.}
\end{itemize}
\item[1.2.] Run script
\item[1.3.] The output file will be saved to the working directory with `{\bf \_pp.mat}' suffix for use in the main DNL script.
\end{itemize}

\noindent {\bf 2. Main program: `DYNAMIC\_NETWORK\_LOADING.m'}
\begin{enumerate}
 \setlength\itemsep{0 em}
\item[2.1.] Set user inputs (at the top of file):
\begin{itemize}
 \setlength\itemsep{0 em}
\item[a)] Location of pre-processed data file (e.g. `.../{\bf Braess\_pp.mat}')
\item[b)] Location of departure rates data file (e.g. `.../{\bf Braess\_dep.mat}'). The path departure vector should be a matrix of dimension $|\mathcal{P}|\times N$ where $|\mathcal{P}|$ is the number of paths and $N$ is the number of time steps.
\item[c)] Number of time steps $N$, which should match the dimension of the path departure vector.
\end{itemize}
\item[2.2.] Run the script `{\bf DYNAMIC\_NETWORK\_LOADING.m}'.
\item[2.3.] The output file will be saved to the working directory with `{\bf \_out.mat}' suffix, containing data required for graphical display.
\end{enumerate}

\noindent {\bf 3. Graphical display: `OPEN$\_$NETWORK$\_$GUI.m'}
\begin{enumerate}
\item[3.1.] Run script (can be run after loading any saved output files into the workspace)
\end{enumerate}

\pagebreak

\subsection*{Appendix 2. Example of the DNL solver}
\begin{figure}[h!]
\centering
\includegraphics[width=.4\textwidth]{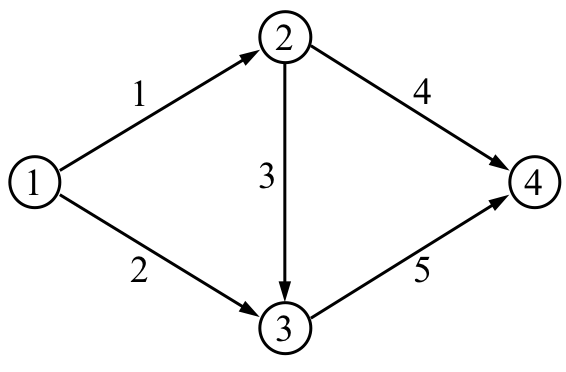}
\caption{The Braess network (example).}
\label{figBraess}
\end{figure}

We consider the Braess network shown in Figure \ref{figBraess} as an illustrative example. This network has four O-D pairs: $(1,\,3)$, $(2,\,3)$, $(1,\,4)$ and $(2,\,4)$, and eight paths:
\begin{itemize}
\item O-D $(1,\,3)$: $p_1=\{1,\,3\}$, $p_2=\{2\}$;
\item O-D $(2,\,3)$: $p_3=\{3\}$;
\item O-D $(1,\,4)$: $p_4=\{1,\,4\}$, $p_5=\{1,\,3,\,5\}$, $p_6=\{2,\,5\}$; and
\item O-D $(2,\,4)$: $p_7=\{4\}$, $p_8=\{3,\,5\}$.
\end{itemize}

\begin{figure}[h!]
\centering
\includegraphics[width=.9\textwidth]{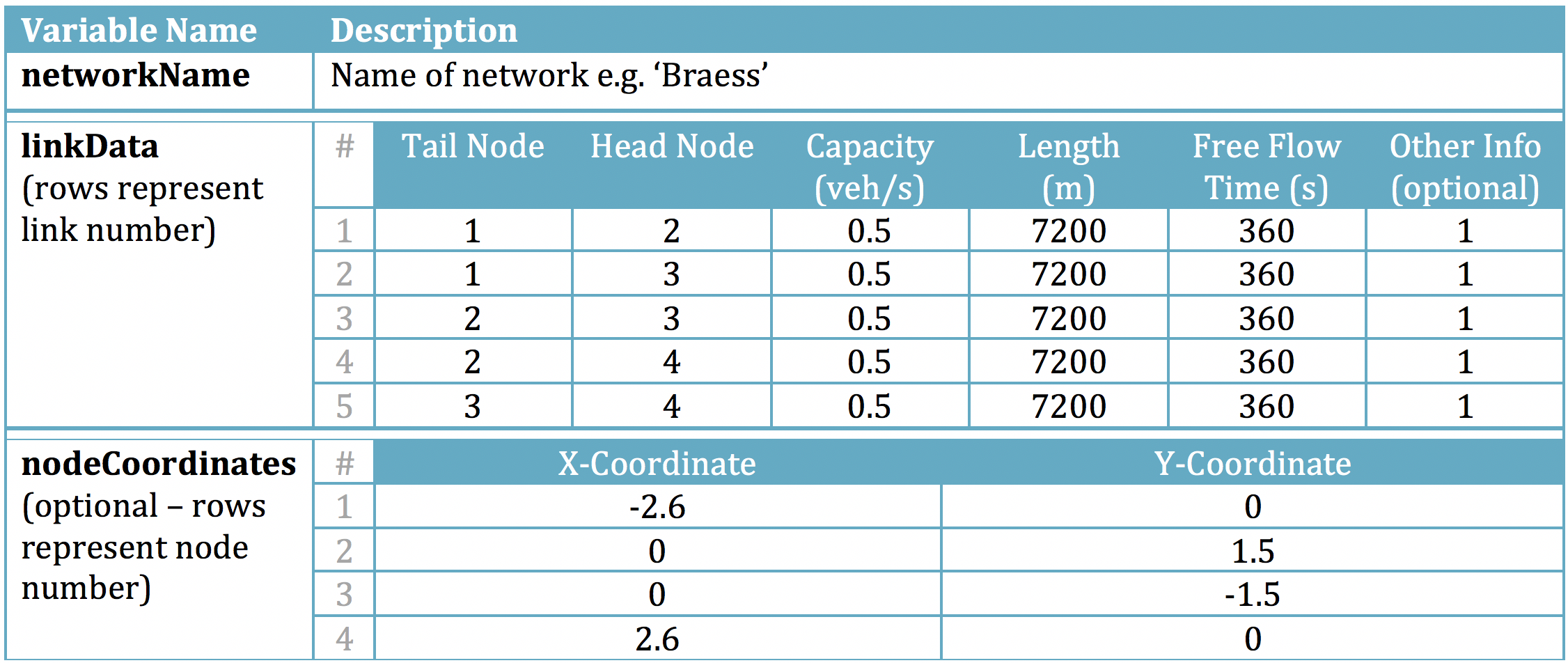}
\caption{Network data for the Braess network (example). Note that, for simplicity, the program assumes that $v=3 w$ where $v$ and $w$ are respectively the speeds of the forward and backward kinematic waves; see \eqref{trifd}. Therefore, the link attributes listed in this figure will identify a unique triangular fundamental diagram.}
\label{figTable1}
\end{figure}

\begin{figure}[h!]
\centering
\includegraphics[width=.9\textwidth]{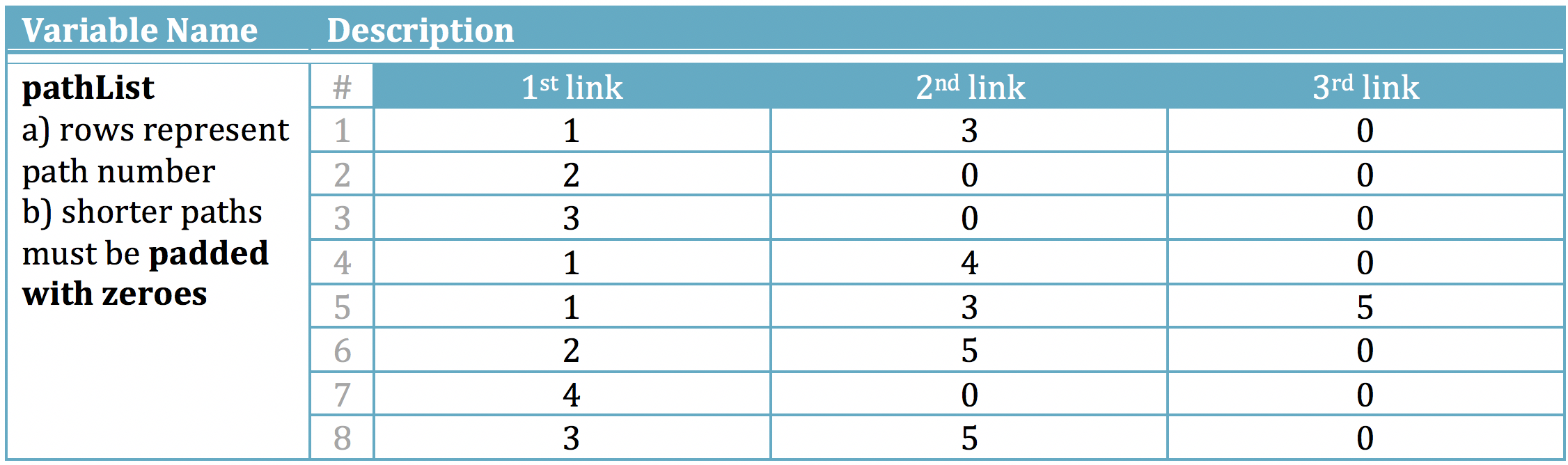}
\caption{Path data for the Braess network (example).}
\label{figTable2}
\end{figure}

\begin{figure}[h!]
\centering
\includegraphics[width=.9\textwidth]{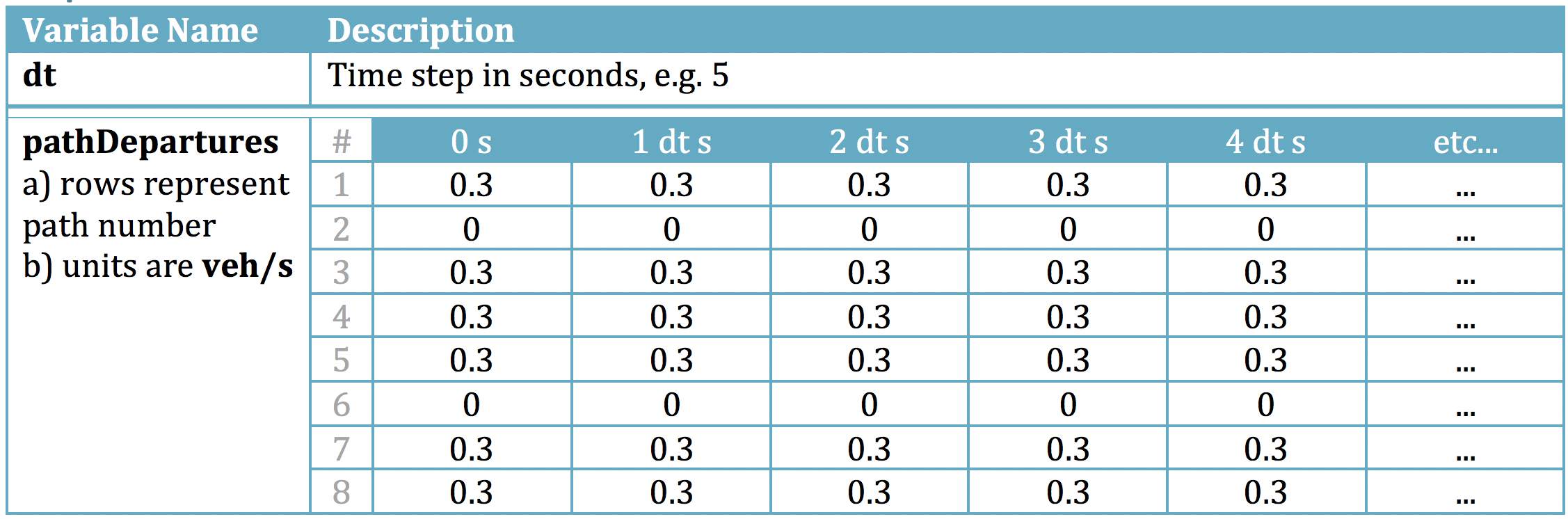}
\caption{Path departure rate vector data for the Braess network (example).}
\label{figTable3}
\end{figure}

\noindent The network data, path data, and path departure vector are illustrated in Figures \ref{figTable1}, \ref{figTable2} and \ref{figTable3}, respectively.

To run the dynamic network loading procedure for the Braess network, we follow the diagram in Figure \ref{figFlow1}: 
\begin{enumerate}
\item Generate two data files, namely `{\bf Braess\_dat.mat}' and `{\bf Braess8\_paths.mat}' according to Figures \ref{figTable1} and \ref{figTable2};\footnote{For large networks, the path set can be generated by using the k-shortest path or Frank-Wolfe algorithms.}
\item Run the script `{\bf CREATE\_NETWORK\_PROPERTIES.m}', which outputs and saves the file `{\bf Braess8\_pp.mat}', which contains all the pre-processed network parameters.
\item Generate the path departure rate vector file `{\bf Braess\_dep.mat}' according to Figure \ref{figTable3}.
\item Run the script `{\bf DYNAMIC\_NETWORK\_LOADING.m}' in the same directory, which then outputs and saves the file `{\bf Braess\_out.mat}'. This file includes path travel times, which is required by the fixed-point algorithm for computing DUEs, and several link variables that are required for the graphical display. 
\item (optional) Run `{\bf OPEN\_NETWORK\_GUI.m}' to open and operate the Graphical User Interface. 
\end{enumerate}

\begin{figure}[h!]
\centering
\includegraphics[width=.98\textwidth]{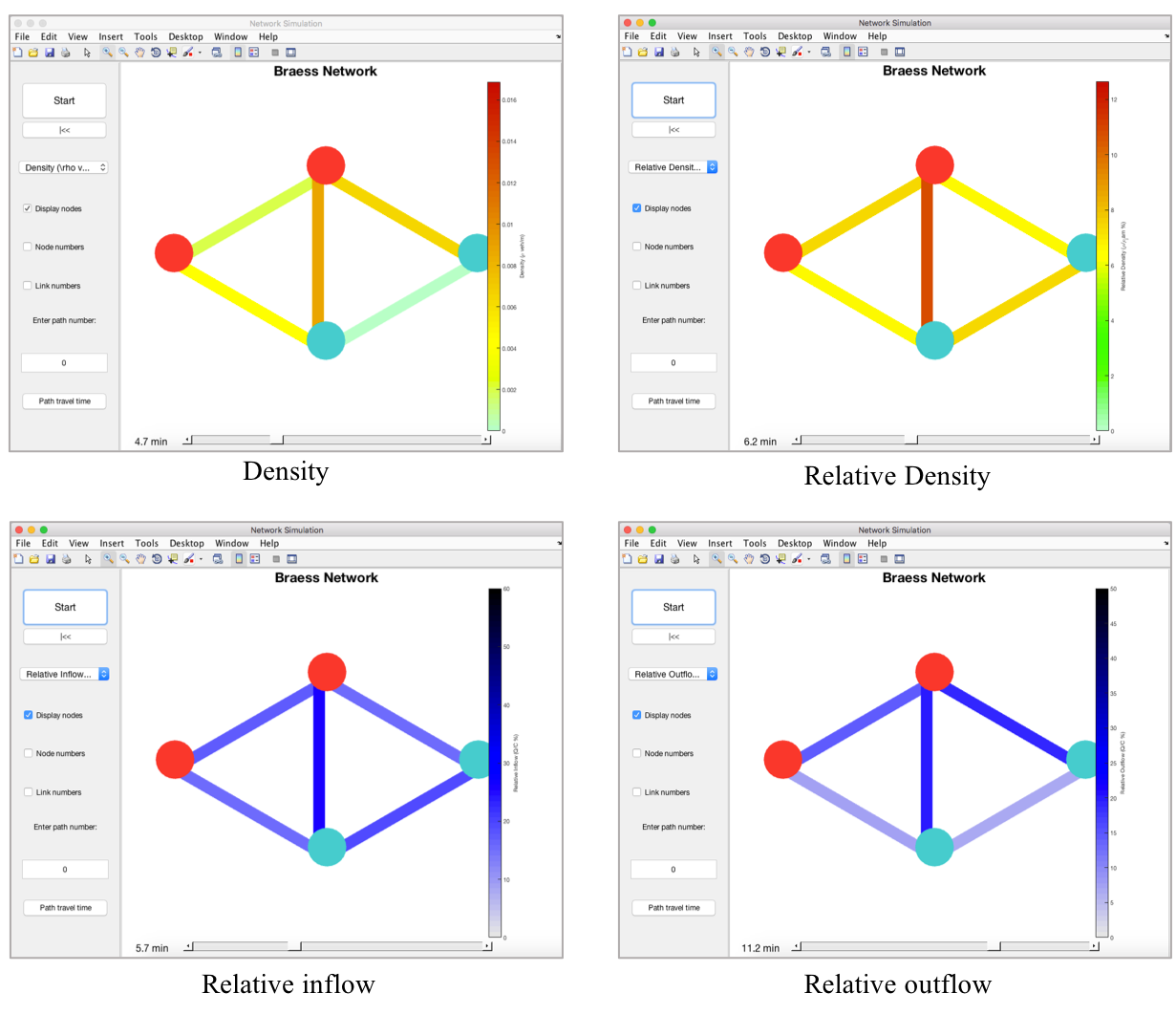}
\caption{Snapshots of the graphical user interface for the Braess network.}
\label{figBraessGUI}
\end{figure}

\noindent Figure \ref{figBraessGUI} illustrates the Graphical User Interface. Four display options are available:
\begin{itemize}
 \setlength\itemsep{0 em}
\item {\it Density} (veh/m): average link density computed as the number of vehicles traveling on the link at any time over the link length.
\item {\it Relative density} (\%): the aforementioned density over the jam density of the link.
\item {\it Relative inflow} (\%): the inflow of the link at any time over the link capacity.
\item {\it Relative outflow} (\%): the outflow of the link at any time over the link capacity.
\end{itemize}

\pagebreak

\subsection*{Appendix 3. Dynamic user equilibrium solver}

\begin{figure}[h!]
\centering
\includegraphics[width=.8\textwidth]{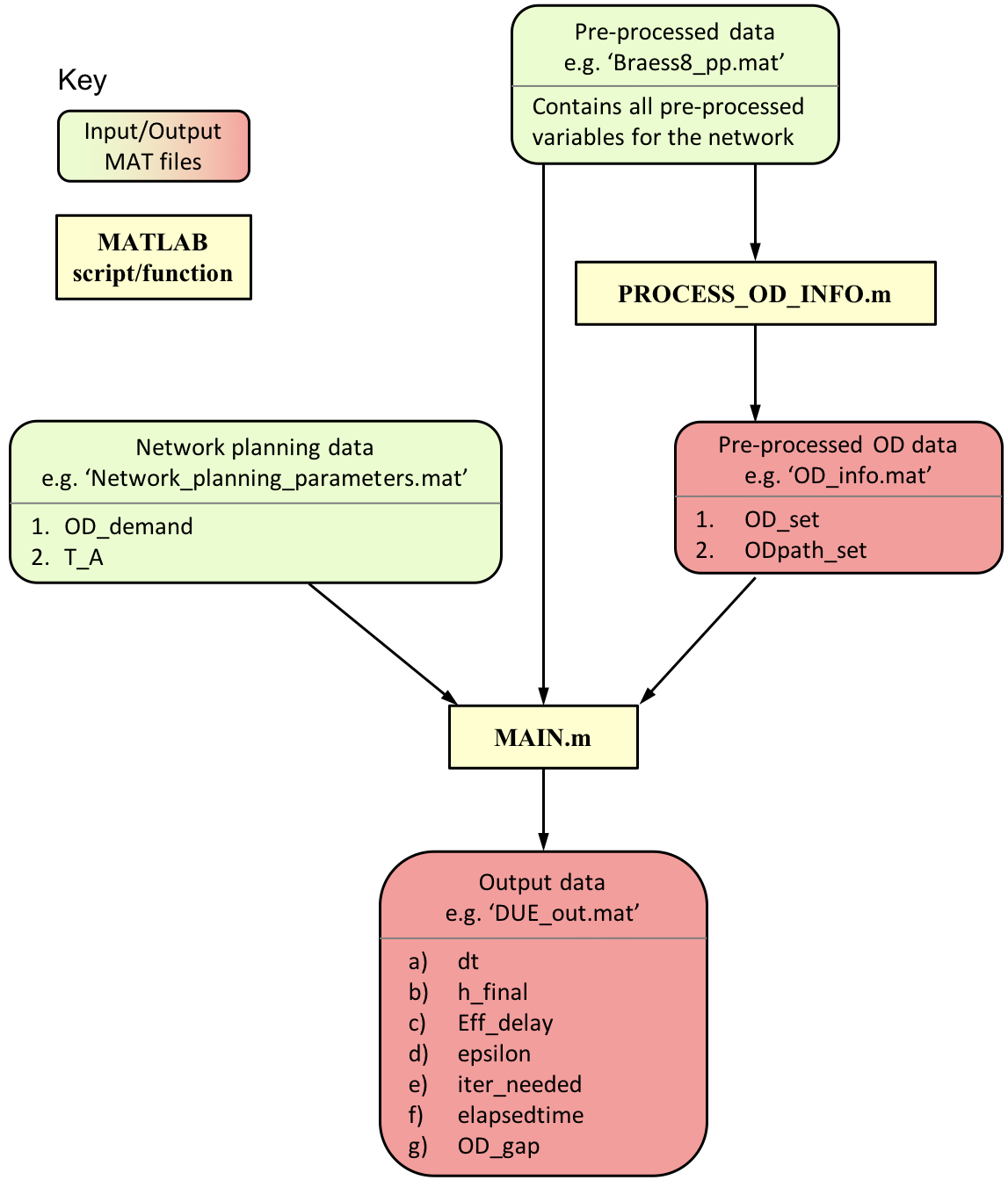}
\caption{Structure and interdependencies of all input and output files within the DUE program. Files in green are user-provided data containing the variables listed.}
\label{figDUEflow}
\end{figure}

The dynamic user equilibrium solver is based on the fixed-point algorithm presented in Section \ref{secFPalgo}. The individual scripts or functions (.m files) and input data files are illustrated in Figure \ref{figDUEflow}.

The MATLAB implementation of the DUE solution procedure consists of three main steps:\\
\noindent {\bf 1. Pre-processing the O-D information: `PROCESS\_OD\_INFO.m'}\\
The purpose of this script is to process the network data (e.g. `{\bf Braess8\_pp.mat}') in order to extract information about all the O-D pairs and their association with the paths. Follow the steps below to run this script.
\begin{itemize}
 \setlength\itemsep{0 em}
\item[1.1.] Prepare the pre-processed network data (e.g. `{\bf Braess8\_pp.mat}') and place it in the same working directory as `{\bf PROCESS\_OD\_INFO.m}'.
\item[1.2.] Run the script `{\bf PROCESS\_OD\_INFO.m}'.
\item[1.3.] The output file, named `{\bf OD\_info.mat}', will be saved to the working directory for use by the DUE solver.
\end{itemize}

\noindent {\bf 2. Specify input file `Network\_planning\_parameters.mat'}\\
This step mainly specifies the O-D demand and target arrival times. The O-D demand data `{\bf OD\_demand}' is a vector of positive numbers, whose dimension is equal to the number of O-D pairs in the network. The target arrival time data `{\bf T\_A}' is a vector of target arrival times, whose dimension is also equal to the number of O-D pairs. This means that each O-D pair is associated with one target arrival time. 

\vspace{0.1 in}

\noindent {\bf 3. Main program: `MAIN.m'}
\begin{itemize}
 \setlength\itemsep{0 em}
\item[3.1.] Set user inputs (at the top of file):
\begin{itemize}
 \setlength\itemsep{0 em}
\item[a)] Time step (in second) of the discretized problem. Note that ideally the time step should be no less than the minimum link free-flow time in the network to ensure numerical stability. However, the code is conditioned to accommodate violation of this rule, in which case a warning message will be displayed. 
\item[b)] Threshold that serves as the convergence criterion of the fixed-point algorithm; see \eqref{chapNum:fptermination}.
\item[c)] Tolerance (indifference band) if bounded rationality (BR) is to be considered \citep{HSF2015}; the default value is set to be 0, which means no BR is considered. Empirically, and understandably, a positive tolerance tends to facilitate convergence of the fixed-point algorithm.
\item[d)] The maximum number of fixed-point iterations, upon which the algorithm will be terminated regardless of convergence.
\item[e)] The step size $\alpha$ of the fixed-point algorithm; see \eqref{eqnFPalgupdate}. Note that $\alpha$ needs to be adjusted for different networks and demand levels, as it has a major impact on convergence with respect to the criterion \eqref{chapNum:fptermination}. Indeed, a very small step size $\alpha$ tends to terminate the fixed-point iterations prematurely while compromising the quality of the approximate DUE solution. On the other hand, a very large $\alpha$ causes undesirable oscillations among the iterates; and it normally takes more iterations to reach convergence. Another consideration, arising from \eqref{eqnFPalgupdate}, is that $\alpha\Psi_p(t,\,h^k)$ should be numerically comparable to $h_p^k(t)$. 
\item[f)] Set the file path for the pre-processed network data (e.g. `.../{\bf Braess8\_pp.mat}').
\end{itemize}

\item[3.2.] Run the script `{\bf MAIN.m}'.

\item[3.3.] The output file, named `{\bf DUE\_out.mat}', will be saved to the working directory. The output includes the following variables:
\begin{itemize}
\item[a)] `{\bf dt}': the discrete time step size.
\item[b)] `{\bf h\_final}': the path flow vector upon convergence or forced termination of the fixed-point algorithm. 
\item[c)] `{\bf Eff\_delay}': the effective path delay vector, whose dimension matches that of `{\bf h\_final}'.
\item[d)] `{\bf epsilon}': the relative gap between two consecutive iterates of the fixed-point algorithm. 
\item[e)] `{\bf iter\_needed}': the number of iterations performed upon termination of the algorithm (either when the convergence criterion is met or the maximum iteration number is reached).
\item[f)] `{\bf elapsedtime}': time taken to run the DUE solver. 
\item[g)] `{\bf OD\_gap}': the travel cost gap for all the O-D pairs; see \eqref{gapdef}.
\end{itemize}
\end{itemize}

\begin{figure}[h!]
\centering
\includegraphics[width=.9\textwidth]{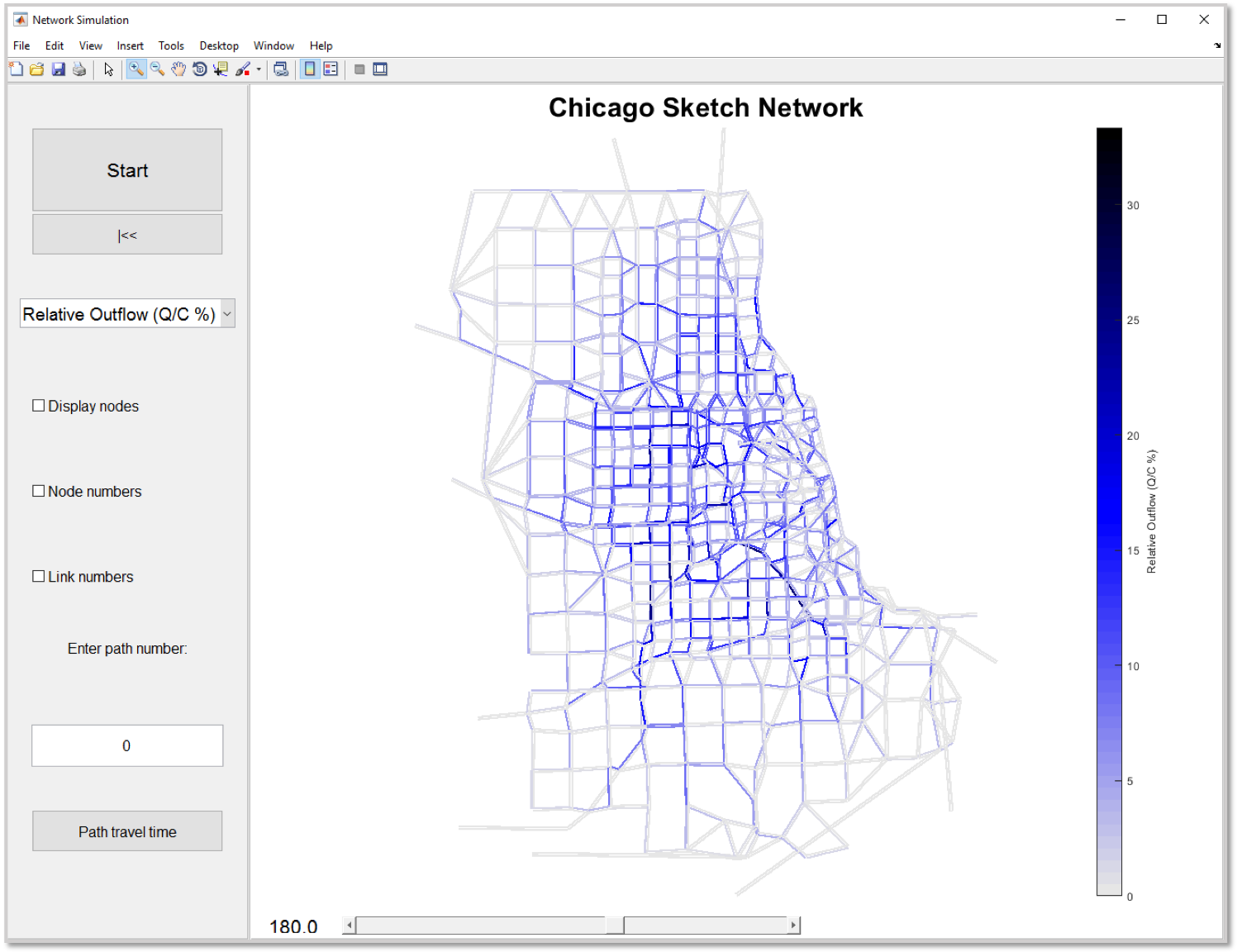}
\caption{Snapshot of the GUI that visualizes the DUE solution on the Chicago Sketch network. The color scale indicates the relative link inflow.}
\label{figSCDUE}
\end{figure}

Figure \ref{figSCDUE} illustrates the DUE solution on the Chicago Sketch network obtained from the MATLAB solver. The codes and data for this network are available at the GitHub repository. 

\pagebreak

\end{document}